\documentclass[10pt,twoside]{article}
\usepackage[latin1]{inputenc}
\usepackage{amsmath}
\usepackage{graphicx}
\usepackage{yhmath}
\usepackage[table]{xcolor}
\usepackage{mathrsfs} 
\usepackage{amssymb}
\usepackage{url}
\usepackage{makecell}

\usepackage{stmaryrd} 

\input epsf
\def\limiten{\renewcommand{\arraystretch}{0.5}
\begin{array}[t]{c}\stackrel{}{\longrightarrow} \\
{\scriptstyle n\rightarrow
\infty}\end{array}\renewcommand{\arraystretch}{1}}

\def\limitepsn{\renewcommand{\arraystretch}{0.5}
\begin{array}[t]{c}\stackrel{a.s.}{\longrightarrow} \\
{\scriptstyle n \rightarrow
\infty}\end{array}\renewcommand{\arraystretch}{1}}

\def\limiteloin{\renewcommand{\arraystretch}{0.5}
\begin{array}[t]{c}\stackrel{{\cal D}}{\longrightarrow} \\
{\scriptstyle n\rightarrow
\infty}\end{array}\renewcommand{\arraystretch}{1}}

\def\limiteproban{\renewcommand{\arraystretch}{0.5}
\begin{array}[t]{c}\stackrel{{\cal P}}{\longrightarrow} \\
{\scriptstyle n\rightarrow
\infty}\end{array}\renewcommand{\arraystretch}{1}}

\numberwithin{equation}{section}

\newtheorem{thm}{Theorem}[section]

\newtheorem{lem}[thm]{Lemma}

\newtheorem{prop}[thm]{Proposition}

\newtheorem{rmrk}[thm]{Remark}

\newcommand{\E}{\ensuremath{\mathbb{E}}}
\newcommand{\R}{\ensuremath{\mathbb{R}}}
\newcommand{\Z}{\ensuremath{\mathbb{Z}}}

\newcommand{\N}{\ensuremath{\mathbb{N}}}

\newcommand{\var}{\ensuremath{\mathrm{Var}}}
\newcommand{\prob}{\ensuremath{\mathbb{P}}}

\definecolor{grisclair}{gray}{0.9}
\font\dsrom=dsrom10 scaled 1200

\renewcommand{\arraystretch}{.8}

\begin{document}
\title{\bf Piecewise autoregression for general integer-valued time series}
 \maketitle \vspace{-1.0cm}
  \begin{center}
   Mamadou Lamine DIOP$^{\text{a, b}}$ \footnote{Supported by the Institute for advanced studies - IAS (Université de Cergy-Pontoise, France),
    the MME-DII center of excellence (ANR-11-LABEX-0023-01) and by the CEA-MITIC (Université Gaston Berger, Sénégal)} and
   William KENGNE$^{\text{a}}$ \footnote{Developed within the ANR BREAKRISK : ANR-17-CE26-0001-01}
 \end{center}

  \begin{center}
 {\it $^{\text{a}}$  THEMA, Université de Cergy-Pontoise, 33 Boulevard du Port, 95011 Cergy-Pontoise Cedex, France \\
 $^{\text{b}}$ LERSTAD, Université Gaston Berger,  BP 234 Saint-Louis, Sénégal\\
  E-mail: mamadou-lamine.diop@u-cergy.fr ; william.kengne@u-cergy.fr\\ }
\end{center}

 \pagestyle{myheadings}
 \markboth{Piecewise autoregression for count time series}{M. L. Diop and W. Kengne}

\textbf{Abstract} :
 This paper proposes a piecewise autoregression for general integer-valued time series. The conditional mean of the process depends on a parameter which is piecewise constant over
 time. We derive an inference procedure based on a penalized contrast that is constructed from the Poisson quasi-maximum likelihood of the model.
 The consistency of the proposed estimator is established.
 From practical applications, we derive a data-driven procedure based on the slope heuristic to calibrate the penalty term of the contrast;
 and the implementation is carried out through the dynamic programming algorithm, which leads to a procedure of $\mathcal{O}(n^2)$ time complexity.
  Some simulation results  are provided, as well as the applications to the US recession data and the number of trades in the stock of Technofirst.

{\em Keywords:} Multiple change-points, model selection,  integer-valued time series, Poisson quasi-maximum likelihood,  penalized quasi-likelihood, slope heuristic.

\section{Introduction}
 We consider a $\N_0$-valued ($\N_0=\N \cup \{0\}$)  process $Y=\{Y_{t},t\in \Z \}$ where the conditional mean
   \begin{equation}\label{model}
     \lambda_t = \lambda_t(\theta^*_t)  = \E(Y_t | \mathcal{F}_{t-1})
   \end{equation}
 is a function (see below) of the whole information $\mathcal{F}_{t-1}$  up to time $t-1$ and of an unknown parameter $\theta^*_t$ belongs to a
  compact subset $\Theta \subset \R^d$ ($d \in \N$).
%
%
%
  The inference in the cases where $\theta^*_t = \theta^*$  is constant or the distribution of $ Y_t| \mathcal{F}_{t-1}$ is known have been studied
  by many authors in several directions; see for instance Fokianos \textit{et al.} (2009), Fokianos and Tjøstheim (2011, 2012),  Douc \textit{et al.} (2017) among others,
  for some recent works.
  We consider here a more general setting where $\theta^*_t$  is piecewise constant and  the distribution of $Y_t| \mathcal{F}_{t-1}$ is unknown.

   \medskip
    \medskip
 We consider  the observations $Y_1, \cdots , Y_n$ generated as in model (\ref{model}) and assume that the parameter $\theta^*_t$ is piecewise constant.
Assume that $\exists K^*  \in \N$, $ \underline{\theta}^* = (\theta^*_1,\cdots,\theta^*_{K^*}) \in \Theta^{K^*}$ and $0 <t^*_1<\cdots < t^*_{K^*-1}<n$
%
%
 such that, $\{Y_{t}, t^*_{j-1} < t \leq t^*_j\}$ is a trajectory of the process  $\{Y_{t,j},t\in \Z\}$ (see Section \ref{sect_PQMLE}) satisfying:
    \begin{equation}\label{model_pwc}
    \E(Y_{t,j} | \mathcal{F}_{t-1}) = f(Y_{t-1,j}, Y_{t-2,j}, \cdots; \theta^*_j), ~  \forall~ t^*_{j-1} < t \leq t^*_j
   \end{equation}
   where     $  \mathcal{F}_{t} = \sigma(Y_{s,j}, s\leq t, ~ j=1,\cdots,K^*-1)  $    is the $\sigma$-field generated by the whole information up to time $t$ and
    $f$ a measurable non-negative function assumed to be known up to the parameter $\theta^*_t$.
 $K^*$ is the number of segments (or regimes) of the model; the $j$th segment corresponds to $\{t^*_{j-1}+1, t^*_{j-1}+2 ,\cdots,t_j^* \}$ and depends on
 the parameter $\theta^*_j$. $t_1^*, \cdots, t_{K^*-1}^*$ are the change-point locations; by convention, $t^*_0 = -\infty$ and $t^*_{K^*} = \infty$.
    To ensure identifiability of the change-point locations, it is reasonable to assume that $\theta^*_j \neq \theta^*_{j+1}$ for $j=1,\cdots,K^*-1$.
    The case $K^*=1$ corresponds to the  model without change.
    In the sequel, we assume that the random variables $Y_t$, $t \in \Z$ have the same (up to the parameter $\theta_t^*$)  distribution $P$ and denote by
    $P\big(\lambda_t(\theta_t^*) \big)$ the distribution of $Y_t |\mathcal{F}_{t-1}$.
 %
 \medskip
 Our main focus of interest is the estimation of the unknown parameters $\big(K^*, (t_j^*)_{1\leq j \leq K^*-1} , (\theta_j^*)_{1\leq j \leq K^*} \big)$
 in the model (\ref{model_pwc}). This can be viewed as a classical model selection problem.
 Assume that the observations $Y_1,\cdots,Y_n$ are generated from (\ref{model_pwc}).
 Let $K_{\max} $ be the upper bound of the number of segments (note that $K_{\max} < n$).
 Denote by  $\mathcal{M}_n$ the set of partitions of $\llbracket 1, n\rrbracket$ into at most $K_{\max}$  contiguous segments.
 Set $m = \{ T_1,\cdots, T_K \}$ a generic element of $K$ segments in $\mathcal{M}_n$.
 %
 Consider the collection $\{\mathcal{S}_{m}, ~ m \in \mathcal{M}_n \}$ where, for a given  $m \in \mathcal{M}_n$, $\mathcal{S}_{m}$ is the families of sequence
  $(\theta_t)$ which are piecewise constant on the partition $m$.
 %
%
 Any $\vartheta =(\theta_t) \in \mathcal{S}_{m}$ depends on the parameter $\underline{\theta} = (\theta_1,\cdots,\theta_K)$ which is the piecewise values of
 $\theta_t$ on each segment.
 Set $\mathcal{S}= \cup_{m \in \mathcal{M}_n}\mathcal{S}_{m}$.
 Denote by $\vartheta$ a generic element of $ \mathcal{S} $, with partition $m$ and parameter $\underline{\theta}$.
 $ | \underline{\theta} | = K$ denote the number of the piecewise segment, also called the dimension of $\vartheta$.
  The true model $\vartheta^*$ with dimension $K^*$, depends on a partition  $m^*$ and the parameter
 $\underline{\theta}^* $.

\medskip
 For any $\vartheta \in \mathcal{S}$, set
 $  \lambda_t^\vartheta = \sum_{k=1}^K  \lambda_t(\theta_k) \textrm{\dsrom{1}}_{t \in T_k}  $
 and denote by  $P(\lambda_t^{\vartheta})$ the distribution of $Y_t|\mathcal{F}_{t-1},\vartheta$ ; let $ p(\cdot | \mathcal{F}_{t-1},\vartheta)= p(\cdot;\lambda_t^\vartheta)$
 be the probability density function of this distribution.
 For $\vartheta \in \mathcal{S}$, let $P_{n,\vartheta}$ be the conditional distribution of $(Y_1,\cdots,Y_n)|\mathcal{F}_{n-1},\vartheta$.
 %
 %
 %
We consider the log-likelihood contrast: $\forall \vartheta \in \mathcal{S} $,
\[
   ~ ~ \gamma_n( \vartheta ) := \gamma_n( P_{n,\vartheta}  ) = - \log P_{n,\vartheta}(Y_1,\cdots,Y_n)  =  - \sum_{t=1}^n \log p(Y_t|\mathcal{F}_{t-1},\vartheta)
             =  -  \sum_{t=1}^n \log p(Y_t;\lambda_t^{\vartheta}) .
  \]
%
%
%
%
%
Thus, the minimal contrast estimator $\widehat{\vartheta}_{m}$ of $\vartheta^*$ on  the collection $\mathcal{S}_{m}$ is obtained by minimizing the contrast $\gamma_n(\vartheta )$
over $\vartheta \in \mathcal{S}_{m}$ ; that is
$ \widehat{\vartheta}_{m} = \underset{ \vartheta \in \mathcal{S}_{m} }{ \text{argmin}} ~\gamma_n(\vartheta)$  .
%
 %
  The main approaches of the model selection procedures take into account the  model complexity and select the estimator
  $\widehat{\vartheta}_{ m_n} $  such that, $m_n$  minimizes  the  penalized criterion
\begin{equation}\label{Criterion-Pen}
\text{crit}_n(m) = \gamma_n(\widehat{\vartheta}_{m}) + \text{pen}_n(m), ~~\text{for~ all} ~ m \in \mathcal{M}_n
\end{equation}
where $pen_n : \mathcal{M}_n \rightarrow \R_+$ is a penalty function, possibly data-dependent.
We now address the following issues.

 \medskip
  (i) \textbf{Semi-parametric setting}.
      Cleynen and Lebarbier (2014 and 2017) recently consider the change-point type problem (\ref{model_pwc}) with i.i.d. observations; in their works, the distribution $P$ is
      assumed to be known and could be Poisson, Negative binomial or belongs to the exponential family distribution.
       From the practical viewpoint, we consider the case where $P$ is unknown and deal with the Poisson quasi-likelihood (see for instance  Ahmad and Francq (2016)).
       So in the sequel, $\gamma_n$ is the Poisson quasi-likelihood contrast and $\widehat{\vartheta}_{m}$ is the Poisson quasi-maximum likelihood estimator (PQMLE).

  \medskip
   (ii) \textbf{Multiple change-point problem  from a non-asymptotic point of view}.
         This question is tacked by model selection approach. Numerous works have been devoted to this issue; see among others, Lebarbier (2005), Arlot and Massart (2009),
         Cleynen and Lebarbier (2014 and 2017),  Arlot and  Celisse (2016).\\
          In this (quasi)log-likelihood framework, it is more usual to consider the Kullback-Leibler risk.
 %
 %
 %
 For any $\vartheta \in \mathcal{S}$, the Kullback-Leibler  divergence between $P_{n,\vartheta^*}$ and $P_{n,\vartheta}$  is
 \begin{align*}
  KL(\vartheta^*,\vartheta) := KL(P_{n,\vartheta^*}, P_{n,\vartheta}) &=  \E \Big[ \log \frac{ P_{n,\vartheta^*}(Y_1,\cdots,Y_n) }{ P_{n,\vartheta}(Y_1,\cdots,Y_n)  } \Big]
   = \sum_{t=1}^n \E \Big[ \log \frac{  p(Y_t|\mathcal{F}_{t-1},\vartheta^*)  }{  p(Y_t|\mathcal{F}_{t-1},\vartheta) } \Big] \\
  & = \sum_{t=1}^n \E \big[ \log p(Y_t;\lambda_t^{\vartheta^*}) \big] - \sum_{t=1}^n \E \big[ \log p(Y_t;\lambda_t^{\vartheta}) \big]
\end{align*}
 where $\E$ denotes the expectation with respect to the true distribution of the observations.
 %
  %
  In the case where $\gamma_n$ is the likelihood contrast, we get $KL(\vartheta^*,\vartheta) = \E[\gamma_n(\vartheta) - \gamma_n(\vartheta^*)]$.
  %
%
%
%
 The "ideal" partition  $m(\vartheta^*)$ (the one whose estimator is closest to $\vartheta^*$ according to the Kullback-Leibler risk) satisfying:
 \[ m(\vartheta^*) = \underset{ m \in \mathcal{M}_n }{ \text{argmin}} ~ \E[ KL(\vartheta^*,\widehat{\vartheta}_{m}) ]  .\]
 The corresponding estimator, $\widehat{\vartheta}_{m(\vartheta^*)}$ called the \emph{oracle}, depends on the true sample distribution, and cannot be computed
 in practice. The goal is to calibrate the penalty term, such that the segmentation $\widehat{m}$ provides an estimator
 $\widehat{\vartheta}_{ \widehat{m}}$
 where the risk of $\widehat{\vartheta}_{\widehat{m}}$ is close as possible to the risk of the \emph{oracle}, namely such that
  \begin{equation}\label{oracle1}
   \E[ KL(\vartheta^*,\widehat{\vartheta}_{\widehat{m}}) ] \leq   C ~\E[ KL(\vartheta^*,\widehat{\vartheta}_{ m(\vartheta^*) }) ]
  \end{equation}
  for a nonnegative constant $C$, expected close to 1. This issue is addressed in the above mentioned papers, and the results obtained are heavily
  relied on the independence of the observations. In our setting here, it seems to be a more difficult task. But, we believe that the coupling method can be used as in
  Lerasle (2011) to overcome this difficulty. We leave this question as the topic of a different research project.

\medskip
(iii) \textbf{Multiple change-point problem  from an asymptotic point of view}.
 The aim here is to consistently estimate the parameters of the change-point model.
 This issue has been addressed by several authors using the classical contrast/criteria optimization or binary/sequential segmentation/estimaion;
 see for instance Bai and Perron (1998), Davis \textit{et al.} (2008),  Harchaoui and Lévy-Leduc (2010), Bardet \textit{et al.} (2012), Davis and Yau (2013), Davis \textit{et al.} (2016),
 Yau and Zhao (2016), Inclan and Tiao (1994), Bai (1997),  Fryzlewicz and Subba Rao (2014), Fryzlewicz (2014), among others, for some advanced towards this issue.
 These works and many other papers in the literature on the asymptotic study of multiple change-point problem are often focussed on continuous valued time series;
 moreover, the case of a large class of semi-parametric model for discrete-valued time series (such as those discussed earlier) have not yet addressed.

 \medskip
 We consider (\ref{model_pwc}) and derive a penalized contrast of type (\ref{Criterion-Pen}). We assume that there exists a partition $\underline{\tau}^*$ of $[0, 1]$ such that $[\underline{\tau}^* n]=m^*$, where $[\underline{\tau}^* n]$ is the corresponding partition of
 $\llbracket 1, n\rrbracket$ obtained from $\underline{\tau}^*$.
 We provide sufficient conditions on the penalty $\text{pen}_n$,  for which  the estimators $\widehat{m}$ and  $\widehat{\vartheta}_{ \widehat{m}}$ are consistent ; that is:
 \[
  \big(|\widehat{m}|, \frac{\widehat{m}}{n} , \widehat{\vartheta}_{\widehat{m}} \big)  \limiteproban
 \big(K^*, \underline{\tau}^* , \vartheta^* \big)
  \]
 where $\frac{\widehat{m}}{n}$ is the corresponding partition of $[0,1]$ obtained from $\widehat{m}$.

  \medskip

  The paper is organized as follows. In Section 2, we set some notations, assumptions and define the Poisson QMLE. In Section 3, we derive the estimation procedure and provide
  the main results. Some simulations results are displayed in Section 4 whereas Section 5 focus on applications on the US recession data and the daily  number of trades in the
  stock of Technofirst. Section 6 provides the proofs of the main results.

 \section{Notations and Poisson QMLE}\label{sect_PQMLE}

%
%

    We set the following classical Lipschitz-type condition  on the function $f$.

    \medskip
    \noindent \textbf{Assumption} \textbf{A}$_i (\Theta)$ ($i=0,1,2$):
    For any $y\in \mathbb{N}_{0}^{\N}$, the function $\theta \mapsto f(y; \theta)$ is $i$ times continuously differentiable on $\Theta$ and
    there exists a sequence of non-negative real numbers  $(\alpha^{(i)}_k)_{k\geq 1} $ satisfying
    $ \sum\limits_{k=1}^{\infty} \alpha^{(0)}_k <1 $ (or $ \sum\limits_{k=1}^{\infty} \alpha^{(i)}_k <\infty $ for $i=1, 2$) ;
   such that for any  $y, y' \in \mathbb{N}_{0}^{\N}$,
  \[ \sup_{\theta \in \Theta  } \left| \frac{\partial^i f(y ; \theta)}{ \theta^i}-\frac{\partial^i f(y' ; \theta)}{\theta^i} \right|
  \leq  \sum\limits_{k=1}^{\infty}\alpha^{(i)}_k |y_k-y'_k|. \]

 \medskip

  In the whole paper, it is assumed that for $j=1,\cdots,K^*$, there exists a stationary and ergodic process $\{Y_{t,j},t\in \Z\}$ satisfying
    \begin{equation}\label{stat_solu}
      \E(Y_{t,j} | \mathcal{F}_{t-1,j}) = f(Y_{t-1,j}, Y_{t-2,j}, \cdots; \theta^*_j), ~  \forall~ t \in \Z
    \end{equation}
   where  $\mathcal{F}_{t,j}=\sigma(Y_{s,j}, s\leq t)$ is the $\sigma$-field generated by $\{ Y_{s,j}, s\leq t\}$; and
    \begin{equation}\label{moment}
    \exists C>0, \epsilon >1, \text{ such that } \forall t \in \Z, ~ ~  \E Y_{t,j}^{1+\epsilon} < C.
   \end{equation}
  $\{Y_{t,j},t\in \Z\}$ is a stationary solution of the $j$th regime. The focus process $Y=\{Y_{t},t\in \Z \}$ is modelled by these stationary regimes ; that is, for any
  $j=1,\cdots,K^*$, $\{Y_{t}, t^*_{j-1} < t \leq t^*_j\}$ is a trajectory of the process  $\{Y_{t,j},t\in \Z\}$.

  \medskip

 \noindent Ahmad and Francq \cite{Francq2016} (Section 3) have discussed about the stationarity and ergodicity issues.
 In many classical integer-valued time series, the assumption \textbf{A}$_0 (\Theta)$ is enough to enable the existence of a stationary and ergodic process
 satisfying (\ref{stat_solu}).

 %

 \subsection{Notations}
   Assume that a trajectory  $(Y_1,\cdots,Y_n)$ of $Y$  is observed; with $0<t^*_1<\cdots < t^*_{K^*-1}<n$.
   By convention $t^*_{0}=-\infty$ and $t^*_{K^*}=\infty$.
   We will use the following notations.
\begin{itemize}
 \item For any finite set $A$, $|A|$ denote the cardinality of $A$.
  \item For $a,b \in \R$ (with $a \leq b$), $\llbracket a, b\rrbracket = \N \cap [a,b]$ is the set of integers between $a$ and $b$.
\item For $K \in \N$,
    $ \mathcal{M}_n(K)=\big\{\underline{t}=(t_1,\ldots, t_{K-1})\; ;\; 0<t_1<\ldots<t_{K-1}<n\big\}$;
    in particular,
    $ \underline{t}^*=\big(t_1^*,\ldots,t_{K^*-1}^*\big)\in \mathcal{M}_n(K^*) $
    is the true vector of the locations of breaks. When $K = 1$, $\mathcal{M}_n(1) $
      corresponds to the model with no break. \\
      In the sequel, any configuration $\underline{t}=(t_1,\ldots, t_{K-1}) \in \mathcal{M}_n(K)$ is also used as a partition  $\{T_1, T_2,\cdots, T_K  \}$ of
     $\llbracket 1, n\rrbracket$ into  $K$  contiguous segments, where $T_1 = \{1, \cdots, t_1 \}$, $T_j = \{t_{j-1}+1, \cdots, t_j \}$ for $j=2,\cdots,K-1$,
     $T_K = \{t_{K-1}+1, \cdots, n \}$.
     In particular, $T_1^* = \{1, \cdots, t_1^* \}$, $T_j^* = \{t_{j-1}^*+1, \cdots, t_j^* \}$ for $j=2,\cdots,K^*-1$,  $T_{K^*} = \{t_{K^*-1}+1, \cdots, n \}$.
      $\mathcal{M}_n(K)$ corresponds to the set of partitions of $\llbracket 1, n\rrbracket$ into  $K$  contiguous segments.
\item For $K \in \N^*$ and  $\underline{t}\in \mathcal{M}_n(K) $ fixed, we set  $n_k=|T_k|$ for $1\leq k \leq K$.
In particular $n_j^*=|T^*_j|$ for $1\leq j \leq K^*$. For $1 \leq k \leq K$  and  $1 \leq j \leq K^*$, let
     $n_{k,j}= |T^*_j \cap T_k|$.
 \item Let $\underline{\theta}^* = (\theta_1^*,\cdots, \theta_{K^*}^*) \in \Theta^{K^*}$ be the vector of the true parameters of the model (\ref{model_pwc}).
 \end{itemize}
 Throughout the sequel, the following norms will be used:
 {\em
\begin{itemize}
 \item $ \|x \|:= \sqrt{\sum_{i=1}^{p} |x_i|^2 } $ for any $x \in \mathbb{R}^{p}$;

 \item  $\left\|f\right\|_{\Theta}:=\sup_{\theta \in \Theta}\left(\left\|f(\theta)\right\|\right)$ for any function $f:\Theta \longrightarrow \mathbb{R}^{d^{\prime}}$;

\item for $x=(x_1,\cdots,x_K) \in \R^K$, $\left\|x\right\|_m= \underset{1 \leq i \leq K }{\max} |x_i|$;

\item  if $Y$ is a random vector with finite $r-$order moments, we set $\left\|Y_t\right\|_r=\E\left(\left\|Y\right\|^r\right)^{1/r}$.
\end{itemize}}


 \subsection{Poisson QMLE}
Let $(Y_{1},\ldots,Y_{n})$ be a trajectory generated from the model (\ref{model_pwc}).
Since the conditional distribution is assumed to be unknown, the likelihood of the model is unknown. The estimation procedure of the parameters $ \theta^*_j$ is based on the Poisson quasi-maximum likelihood introduced by Ahmad and Francq (2016).
The conditional Poisson  (quasi)log-likelihood of the model (\ref{model_pwc}) computed on a segment $T \subset\{1,\dots,n\}$ is given  (up to a constant) by
%
%
\begin{equation}\label{logvm}
\widehat{L}_n(T, \theta) := \sum_{t \in T}(Y_t\log \widehat{\lambda}_t(\theta)- \widehat{\lambda}_t(\theta)) = \sum_{t\in T} \widehat{\ell}_t(\theta) ~ \text{ with }
 \widehat{\ell}_t(\theta) =  Y_t\log \widehat{\lambda}_t(\theta)- \widehat{\lambda}_t(\theta)
 \end{equation}
  where $ \widehat{\lambda}_t(\theta) =  \widehat{f}^\theta_{t}= f(Y_{t-1}, \cdots Y_{1},0,\cdots,0; \theta)$.\\
 According to (\ref{logvm}), the Poisson quasi-likelihood estimator (PQMLE) of $ \theta^*_j$ computed on $T$ is defined by
 \begin{equation}\label{emv}
  \widehat{\theta}_n(T) :=  \underset{\theta\in \Theta}{\text{argmax}} (\widehat{L}_n(T,\theta)).
  \end{equation}

  \medskip

  \noindent Now, for $j = 1,\cdots,K^*$, define the Poisson  (quasi)log-likelihood of the $j$th regime by
\[ L_{n,j}(T^*_j, \theta) := \sum_{t \in T^*_j}(Y_{t,j}\log \lambda_{t,j}(\theta)- \lambda_{t,j}(\theta)) = \sum_{t\in T^*_j} \ell_{t,j}(\theta) ~ \text{ with }
  \ell_{t,j}(\theta) = Y_{t,j}\log \lambda_{t,j}(\theta)- \lambda_{t,j}(\theta)\]
  where $ \lambda_{t,j}(\theta) =f^\theta_{t,j}= f(Y_{t-1,j}, Y_{t-2,j}, \cdots ; \theta)$.
  It can be approximated by
\begin{equation}\label{logvm_j}
 \widehat{L}_{n,j}(T^*_j, \theta) := \sum_{t \in T^*_j}(Y_{t,j}\log \widehat{\lambda}_{t,j}(\theta)- \widehat{\lambda}_{t,j}(\theta)) = \sum_{t\in T^*_j} \widehat{\ell}_{t,j}(\theta) ~ \text{ with }
  \widehat{\ell}_{t,j}(\theta) = Y_{t,j}\log \widehat{\lambda}_{t,j}(\theta)- \widehat{\lambda}_{t,j}(\theta)
  \end{equation}
  where $ \widehat{\lambda}_{t,j}(\theta) = \widehat{f}^\theta_{t,j}=f(Y_{t-1,j} , \cdots, Y_{t^*_{j -1} +1,j},0 \cdots 0 ; \theta)$.\\
  According to (\ref{logvm_j}), the PQMLE  of $ \theta^*_j$ computed on $T^*_j$ is defined by
 \begin{equation}\label{emv_j}
  \widetilde{\theta}_n(T^*_j) :=  \underset{\theta\in \Theta}{\text{argmax}} (\widehat{L}_{n,j}(T^*_j,\theta)).
  \end{equation}
%
%
 %
 %
  To avoid problems of parameter identifiability and to study asymptotic normality of the PQMLE, we shall assume: \\
 (\textbf{A0}): for all  $(\theta, \theta')\in \Theta^2$,
 $ \Big( f(Y_{t-1}, Y_{t-2}, \cdots; \theta)= f(Y_{t-1}, Y_{t-2}, \cdots; \theta')  \ \text{a.s.} ~ \text{ for some } t \in \N \Big) \Rightarrow ~ \theta = \theta'$ ;
 moreover, $\exists  \underline{c}>0$ such that $\displaystyle \inf_{ \theta \in \Theta} f(y_{1}, y_{2}, \cdots; \theta)  \geq \underline{c}$, for all $ y \in  \N_0^{\N} $.

 \medskip

\noindent In order to ensure the consistency and asymptotic normality of the PQMLE, we set the following assumptions for each segment $j= 1,\cdots,K^*$
(see also \cite{Francq2016}):
\begin{enumerate}
    \item [(\textbf{A1}):] $\theta^*_j$ is an interior point of $\Theta \subset \mathbb{R}^{d}$;
  \item [(\textbf{A2}):]  $a_{t,j} \longrightarrow 0$ and $Y_{t,j} a_{t,j} \longrightarrow 0$  as $t\rightarrow \infty$, where $a_{t,j}=\underset{\theta \in \Theta }{\sup} \left| \widehat{\lambda}_{t,j}(\theta) -\lambda_{t,j}(\theta)\right|$;
     \item [(\textbf{A3}):] $J_j=\E \Big[ \frac{1}{\lambda_{t,j}(\theta^*_j)}  \frac{\partial \lambda_{t,j}(\theta^*_j)}{ \partial \theta} \frac{\partial \lambda_{t,j}(\theta^*_j)}{ \partial \theta'}  \Big] <\infty$
           ~and~ $I_j=\E \Big[ \frac{\var(Y_{t,j}|\mathcal{F}_{t-1})}{\lambda^2_{t,j}(\theta^*_j)}  \frac{\partial \lambda_{t,j}(\theta^*_j)}{ \partial \theta} \frac{\partial \lambda_{t,j}(\theta^*_j)}{ \partial \theta'}  \Big] <\infty$;
     \item [(\textbf{A4}):] for all $c' \in \R$, $c' \frac{\partial \lambda_{t,j} (\theta^*_j)}{\partial    \theta}=0$ a.s   $\Rightarrow ~ c'=0$;
     \item [(\textbf{A5}):] there exists a neighborhood $V(\theta^*_j)$ of $\theta^*_j$ such that: for all $i, k \in \left\{1,\cdots,d\right\} $, \[\E \left[ \sup_{\theta \in V(\theta^*_j)} \left|  \frac{\partial^2 }{ \partial \theta_i \partial \theta_k } \ell_{t,j}(\theta) \right|\right]<\infty ; \]
  %
    \item [(\textbf{A6}):]   $b_{t,j}$, $b_{t,j} Y_{t,j}$ and $a_{t,j}d_{t,j} Y_{t,j}$ are of order $O(t^{-h})$ for some $h>1/2$, where
    \[b_{t,j}=\underset{\theta \in \Theta }{\sup} \left\{\E \left[\left\|\frac{\partial \widehat{\lambda}_{t,j} (\theta)}{ \partial \theta} -\frac{\partial \lambda_{t,j}(\theta)}{ \partial \theta}  \right\| \right]\right\}~~\textrm{and}~~d_{t,j}=\underset{\theta \in \Theta }{\sup}\max \left\{ \E\left[\left\| \frac{1}{\widehat{\lambda}_{t,j}(\theta)} \frac{\partial \widehat{\lambda}_{t,j} (\theta)}{ \partial \theta} \right\|\right],  \E\left[ \left\| \frac{1}{\lambda_{t,j}(\theta)} \frac{\partial \lambda_{t,j} (\theta)}{ \partial \theta} \right\|\right] \right\}.  \]
\end{enumerate}
 These aforementioned assumptions hold for many classical models, see Ahmad and Francq \cite{Francq2016}. These authors have established that the estimator  $\widetilde{\theta}_n(T^*_j)$ is strongly consistent, for each regime $j \in \{1,\cdots,K^* \}$ ; that is,
\[ \widetilde{\theta}_n(T^*_j) \limitepsn  \theta^*_j . \]
They have also proved the asymptotic normality of $\widetilde{\theta}_n(T^*_j)$ ; that is,
 \[\sqrt{n^*_j}(\widetilde{\theta}_n(T^*_j)-\theta^*_j) \limiteloin \mathcal{N}(0,\Sigma_j), ~\forall j =1,\cdots,K^*, \]
 where $\Sigma_j:=J^{-1}_j I_j J^{-1}_j$.
 Under the above assumptions, for any $j=1,\cdots ,K^*$, the matrix $\Sigma_j$ can be consistently estimated by (see \cite{Francq2016})
 \begin{align}\label{eq_Sigma_j}
 &\widehat \Sigma_j = \widehat J^{-1}_j \widehat I_j \widehat J^{-1}_j, \text{ where} \\
  \nonumber &\widehat J_j =  \frac{1}{n}\sum_{t=1}^{n}\frac{1}{\widehat \lambda_{t,j}(\widetilde{\theta}_n(T^*_j))}  \frac{\partial \widehat \lambda_{t,j}(\widetilde{\theta}_n(T^*_j))}{ \partial \theta} \frac{\partial \widehat \lambda_{t,j}(\widetilde{\theta}_n(T^*_j))}{ \partial \theta'} ,\\
  \nonumber &\widehat I_j = \frac{1}{n}\sum_{t=1}^{n} \Big(\frac{Y_t}{\widehat \lambda_{t,j}(\widetilde{\theta}_n(T^*_j))} -1 \Big)^2 \,  \frac{\partial \widehat \lambda_{t,j}(\widetilde{\theta}_n(T^*_j))}{ \partial \theta} \frac{\partial \widehat \lambda_{t,j}(\widetilde{\theta}_n(T^*_j))}{ \partial \theta'} .
 \end{align}

\noindent If we consider the process $\{Y_{t},t\in \Z \}$, these properties are also verified on the segment $T^*_1$ since it is  easy to see that
$\{ (Y_t, \lambda_t), t\in T^*_1 \}$ is a stationary process while $\{ (Y_t, \lambda_t), t > t^*_1 \}$ is not.


 \medskip

\noindent The following proposition establishes the consistency of the estimator $\widehat{\theta}_n(T^*_j)$,  for any $j \in \{ 1,\cdots K^* \}$.
\begin{prop}\label{prop1}
Assume that (\textbf{A0})-(\textbf{A2}) and (\textbf{A$_0 (\Theta)$}) hold. Then
\[
 \widehat{\theta}_n(T^*_j) \limitepsn \theta^*_j , ~\forall j=1,\cdots,K^*.
\]
\end{prop}
The results of this Proposition have been obtained by Ahmad and Francq (2016) when $(Y_t, \lambda_t)$  is strictly stationary.

\section{Estimation procedure and main results}
 In this section, we carry out the estimation of the number of breaks $K^*-1$ and the instants of breaks $\underline{t}^*$ by using a penalized contrast. Some asymptotic studies are also reported.
 %
\subsection{Penalized Poisson quasi-likelihood estimator}
For any configuration of periods $K \geq 1$, $\underline{t} \in \mathcal{M}_n(K)$ and $ \underline{\theta}= (\theta_1,\cdots, \theta_{K}) \in \Theta^{K^*}$, we define the contrast
 \begin{equation}\label{Jn}
  (QLIK)~~~~~~\widehat{J}_n(K,\underline{t},\underline{\theta}):= -2\sum_{k=1}^{K}  \widehat L_n(T_k,\theta_k).
 \end{equation}
According to the proprieties of the PQMLE (see \cite{Francq2016}), when $K^*$ is known, a natural estimator of
$(\underline{t}^*, \underline{\theta}^*)=((t^*_j)_{1\leq j \leq K^*-1}, (\theta^*_j)_{1\leq j \leq K^*})$ for the model (\ref{model_pwc})
is therefore the PQMLE on every interval $[t_j+1,\cdots,t_{j+1}]$ and every parameters $\theta_j$ for $1\leq j \leq K^*$.
But, since $K^*$ is assumed to be unknown, we cannot directly use such method. To take into account the estimation of $K^*$,
the most classical solution is to penalize the contrast by an additional term $\kappa_n K$, where $\kappa_n$ represents a regularization parameter.\\
 Now, define the penalized contrast $QLIK$, called $penQLIK$, by

 \begin{equation}\label{penQLIK}
 (penQLIK)~~~~~~\widetilde{J}_n(K,\underline{t},\underline{\theta}):= \widehat{J}_n(K,\underline{t},\underline{\theta}) + \kappa_n K,
 \end{equation}
 with $ \kappa_n \leq n$ and $\kappa_n \limiten+\infty$. \\
\noindent
The estimator of $\left(K^* ,\underline{t}^*,\underline{\theta}^*  \right)$ is defined as one of the minimizers of the penalized contrast:
\begin{equation}\label{em_Kt_Theta}
  \left(\widehat{K}_n ,\widehat{\underline{t}}_n,\widehat{\underline{\theta}}_n  \right)
     \in   \underset{ 1 \leq K \leq  K_{\max}}{\text{argmin}}~~  \underset{  (\underline{t}, \underline{\theta})    \in \mathcal{M}_n(K) \times \Theta^{K}}{\text{argmin}} \left( \widetilde {J}_n(K,\underline{t},\underline{\theta})\right) ~~\textrm{and}~~ \widehat{\underline{\tau}}_n =\frac{\widehat{\underline{t}}_n}{n}.
  \end{equation}
~~\\
We will adapt the slope heuristic procedure to calibrate the penalty term from data (see  Baudry \textit{et al.} (2010)). In this procedure, the criteria $QLIK$  is a linear transformation of the penalty
(here the number of periods $K$) for the most complex models (with $K$ close to $K_{\max}$). This slope should be close to $-\kappa_n/2$.
The slope estimation procedure considers only the linear part of $-QLIK(K)$ with $1 \leq K \leq K_{\max}$.
Note that, in practice, a  numerical algorithm can be used to compute the estimator on each segment; therefore, a minimum size is needed for the numerical computation of the criteria.
Thus, we consider only the periods of length larger than some $u_n$ and we can a priori fix $K_{\max}$ smaller than $[n/u_n]$.
  The complete procedure can be summarized as follows:
\begin{itemize}
    \item [1.] For each $1 \leq K \leq  K_{\max}$, draw $\left( K, -\min_{\underline{t}, \underline{\theta}} QLIK(K)\right)$. Then
compute the slope of the linear part: this slope is $\widehat{\kappa}_n/2$.
\item [2.] Using $\kappa_n=\widehat{\kappa}_n$, draw $\left( K, -\min_{\underline{t}, \underline{\theta}} penQLIK(K)\right)_{1 \leq K \leq  K_{\max}}$. This curve has a global minimum at $\widehat{K}_n$.
\end{itemize}

\subsection{Asymptotic behavior}
 Under some assumptions, we will establish the asymptotic behavior of the estimator  $\left(\widehat{K}_n,\widehat{\underline{t}}_n,\widehat{\underline{\theta}}_n  \right)$.
Throughout this article, we set the following classical assumptions in the problem of break detection:

\medskip
\noindent \textbf{Assumption B.} $\underset{ 1 \leq j \leq  K^*-1}{\min} \left\| \theta^*_{j+1}-\theta^*_j\right\|>0$. Also,
 there exists a vector $\underline{\tau}^* =(\tau^*_1,\cdots,\tau^*_{K-1})$ with  $0<\tau^*_1<\cdots<\tau^*_{K-1}<1$, called the vector of breaks such that
  $t^*_j=\left[n\tau^*_j\right]$, for $j=1,\cdots,K$   (where $\left[\cdot\right]$ is the integer part).

\medskip

The following theorem gives the consistency of the estimator $\left(\widehat{K}_n ,\widehat{\underline{t}}_n,\widehat{\underline{\theta}}_n  \right)$.
\begin{thm}\label{th1}
Assume $K_{\max}>K^*$ and (\textbf{A0})-(\textbf{A2}), \textbf{B}.
If \textbf{A$_0 (\Theta)$} holds and $(\kappa_{\ell})$ satisfies
\begin{equation}\label{eq_th1}
 \sum_{\ell \geq 1} \frac{1}{ \kappa_{\ell} } \sum_{k \geq \ell } \alpha_k^{(0)}   < \infty,
 \end{equation}
then
\[  \left(\widehat{K}_n ,\widehat{\underline{\tau}}_n,\widehat{\underline{\theta}}_n  \right)  \limiteproban  \left(K^*,\underline{\tau}^*,\underline{\theta} ^*\right).\]
\end{thm}
 %

 \medskip
By convention, throughout the sequel, if the vectors $\widehat{\underline{t}}_n$  and $\underline{t}^*$  do not have the same length, complete the shorter of the two vectors with $0$ before computing the norm $\left\| \widehat{\underline{t}}_n -\underline{t}^*\right\|_m$.
 The following theorem establishes the rates of convergence of the estimators  $\widehat{\underline{\tau}}_n$.
\begin{thm}\label{th2}
Assume $K_{\max}>K^*$ and (\textbf{A0})-(\textbf{A2}), \textbf{B}.
If \textbf{A$_i (\Theta)$} ($i=0,1,2$), (\ref{moment}) (with $\epsilon > 2$) hold and $(\kappa_{\ell})$ satisfies
\begin{equation}\label{eq_th2}
 \sum_{\ell \geq 1} \frac{1}{ \kappa_{\ell} } \sum_{k \geq \ell } \alpha_k^{(i)}   < \infty,
 \end{equation}
then the sequence $\left( \left\| \widehat{\underline{t}}_n -\underline{t}^*\right\|_m\right)_{n>1}$ is uniformly tight in probability, that is,

\[ \underset{\delta \rightarrow \infty }{\lim}~ \underset{n \rightarrow \infty }{\lim}  \mathbb{P} \left(  \left\| \widehat{\underline{t}}_n -\underline{t}^*\right\|_m >  \delta\right)=0.
\]
\end{thm}

 \medskip
Now, we give the convergence in distribution of the estimator of $\widehat{\underline{\theta}}_n$. By convention, if $\widehat{K}_n<K^*$, set $\widehat{T}_j= \widehat{T}_{\widehat{K}_n}$, for $j \in \left\{\widehat{K}_n,\cdots,K^* \right\}$.
The following theorem establishes the asymptotic normality of $\widehat{\theta}_n(\widehat{T}_j)$.

\begin{thm}\label{th3}
 Assume  $K_{\max}>K^*$  and (\textbf{A0})-(\textbf{A6}) and \textbf{B}.
 If (\ref{moment}) (with $\epsilon > 2$), \textbf{A$_i (\Theta)$} ($i=0,1,2$) hold, such that $(\kappa_{\ell})$ satisfies
 \begin{equation}\label{eq_th3}
 \sum_{\ell \geq 1} \max \big(\frac{1}{ \kappa_{\ell} }, \frac{1}{ \sqrt{\ell}}\big) \sum_{k \geq \ell } \alpha_k^{(i)}   < \infty,
 \end{equation}
 %
  then
 \[
 \sqrt{n^*_j}\left(\widehat{\theta}_n(\widehat{T}_j)-\theta^*_j\right) \limiteloin \mathcal{N}_d \left(0,\Sigma_j\right),~~ \forall\, j=1,\cdots,K^*,
 \]
 where $\Sigma_j:=J^{-1}_j(\theta^*_j) I_j (\theta^*_j) J^{-1}_j(\theta^*_j)$ with
\[
J_j(\theta^*_j)=\E \Big[ \frac{1}{\lambda_{t,j}(\theta^*_j)}  \frac{\partial \lambda_{t,j}(\theta^*_j)}{ \partial \theta} \frac{\partial \lambda_{t,j}(\theta^*_j)}{ \partial \theta'}  \Big]
           ~~ \textrm{and}~~  I_j(\theta^*_j)=\E \Big[ \frac{\var(Y_{t,j}|\mathcal{F}_{t-1})}{\lambda^2_{t,j}(\theta^*_j)}  \frac{\partial \lambda_{t,j}(\theta^*_j)}{ \partial \theta} \frac{\partial \lambda_{t,j}(\theta^*_j)}{ \partial \theta'}  \Big].
 \]
\end{thm}
\medskip

\begin{rmrk}\label{rmk1}
 The conditions on the regularization parameters $(\kappa_n)_{n\in \N}$ can be obtained if the Lipschitzian coefficients of $f(\cdot \,;\, \theta)$ and its
 derivatives are bounded by a geometric or Riemanian sequence:
 %
\begin{enumerate}
    \item  the geometric case: if $\alpha^{(i)}_k= O(a^k)$ ($i=0,1,2$) with $0 \leq a <1$, then any choice of  $(\kappa_n)_{n\in \N}$ such that $\kappa_n \leq n$
    and $\kappa_n \rightarrow \infty$ satisfies (\ref{eq_th1}) and (\ref{eq_th3}) (for instance $\kappa_n$ of order $\log n$ as in the BIC approach).
 \item the Riemanian case:  if $\alpha^{(i)}_k= O(k^{-\gamma})$ ($i=0,1,2$) with $\gamma >3/2$,
\begin{itemize}
    \item if $\gamma >2$, then the conditions (\ref{eq_th1}) and (\ref{eq_th3}) hold for any choice of $(\kappa_n)_{n\in \N}$ such that $\kappa_n \leq n$ and $\kappa_n \rightarrow \infty$.
    \item if $3/2<\gamma \leq 2$, then one can choose any sequence such that  $\kappa_n=O(n^\delta)$ with $\delta >2-\gamma$ or $\kappa_n=O(n^{2-\gamma} (\log n)^\delta)$ with $\delta >1$.
\end{itemize}
\end{enumerate}
\end{rmrk}


\section{Some simulations results}
In this section, we implement the procedure on the R software (developed by the CRAN project). We will restrict our attention to the estimation of
the vector $(K^*,\underline{t}^*)$ ; {\it i.e} the number of segments $K^*$ and the instants of breaks $\underline{t}^*$.
For the performances of the estimator of the parameter $\underline{\theta}^*$, we refer to the works of Ahmad and Francq (2016).
For each process, we generate $100$ replications following the scenarios considered. The estimated number of segments is computed by
using $QLIK$ criteria penalized with $\kappa_n= \widehat{\kappa}_n$, $\kappa_n= \log n$ and $\kappa_n= \log n^{1/3}$. The value of the estimator $\widehat{\kappa}_n$
is calibrated by using the slope estimation procedure (see Baudry \textit{et al.} (2010)) as  described above. Once the regularization parameter $\kappa_n$ obtained, the dynamic programming
algorithm is used to minimize the criteria.
With this algorithm, the complexity of the procedure declines from $\mathcal{O}(n^{K_{\max}})$ to $\mathcal{O}(n^2)$.

\subsection{Implementation procedure}
We give the steps of the dynamic programming algorithm for computing the number of segments $\widehat{K}_n$ and the optimal configuration of the breaks $\widehat{\underline{t}}_n$.
This algorithm is such that if $(t_1, \cdots , t_{K-1}, t)$ represents the optimal configuration of $Y_1, \cdots ,Y_t$ into $K$ segments, then $(t_1, \cdots , t_{K-1})$ is the optimal configuration of $Y_1, \cdots ,Y_{t_{K-1}}$ into $K -1$ segments.
 Assume that the regularization parameter $\kappa_n$ is known and let $ML$ be the upper triangular matrix of dimension $n \times n$ with $ML_{i,l}= \widehat{L}(T_{i,l}, \widehat{\theta}_n(T_{i,l}))$, where $T_{i,l}=\left\{ i,i+1,\cdots,l\right\}$,  for $1 \leq i \leq l \leq n$.
 We summarize the implementation of the procedure as follows:
\begin{itemize}
    \item  \underline{\em The number of segments $\widehat{K}_n$}:
     Let $C$ be an upper triangular matrix of dimension $K_{\max} \times n$. For $1 \leq K \leq  K_{\max}$ and $K \leq t \leq n$, $C_{K,t}$ will be the minimum penalized criteria of $Y_1, \cdots , Y_t$ into $K$ segments. For $t = 1, \cdots, n$,  $C_{1,t} = -2 ML_{1,t}+\kappa_n$ and the relation $C_{K+1,t} = \min_{K \leq l \leq t-1} \left(C_{K,l} -2 ML_{l+1,t} + \kappa_n \right)$ is satisfied. Hence, $\widehat{K}_n = \text{argmin}_{1 \leq K \leq K_{\max}} \left(C_{K,n}\right)$.

     \item  \underline{\em The change-point locations $\widehat{\underline{t}}_n$}:
      Let $Z$ be an upper triangular matrix of dimension $(K_{\max}-1) \times n$. For $1 \leq  K \leq (K_{\max}-1)$ and $K+1 \leq  t \leq n$, $Z_{K,t}$ will be the $K$th potential break-point of $Y_1, \cdots , Y_t$. Therefore, the relation  $Z_{K,t} = \min_{K \leq l \leq t-1} \left(C_{K,l} -2 ML_{l+1,t} + \kappa_n \right)$ is satisfied for $K = 1, \cdots  ,(K_{\max}-1)$.  The break-points are obtained as follow: set $\widehat{t}_{\widehat{K}_n}=n$ and for $K = \widehat{K}_n - 1, \cdots , 1$,   $\widehat{t}_{K}=Z_{K,\widehat{t}_{K+1}}$.
\end{itemize}

\subsection{Results of simulations}


\subsubsection{Poisson-INARCH models}
We consider the problem (\ref{model_pwc}) for a Poisson-INARCH(1), {\it i.e.}  $(Y_1,\cdots,Y_n)$ is a trajectory of the process  $Y=\{Y_{t},t\in \Z \}$ satisfying:
\begin{equation}\label{Poisson_pwc_INARCH}
    Y_{t}|\mathcal{F}_{t-1} \sim \mathcal{P}(\lambda_{t})~~;~~  \lambda_{t}= 
    f(Y_{t-1}, Y_{t-2}, \cdots; \theta^*_j)=\alpha^{(j)}_0 + \alpha^{(j)} Y_{t-1} , ~~  ~  \forall~ t \in T^*_j,~~   \forall~ j \in \left\{1,\cdots,K^*\right\}.
   \end{equation}
The parameter vector  is $\theta^*_j=(\alpha^{(j)}_0 ,\alpha^{(j)})$, for all $j \in \{1,\cdots,K^*\}$.\\

\noindent For $n=500$ and $n=1000$, we generate a sample $(Y_1, \cdots ,Y_n)$ in the following situations:

\begin{itemize}
    \item {\bf scenario IA}$_0$: $\theta^{*}_1=(0.5,0.6)$ is constant ($K^*=1$) ;
    \item {\bf scenario IA}$_1$: $\theta^*_1=(0.5,0.6)$ changes to $\theta^*_2=(1.0, 0.6)$ at $t^*=0.5n$ ($K^*=2$) ;
    \item {\bf scenario IA}$_2$: $\theta^*_1=(0.5,0.6)$ changes to $\theta^*_2=(1.0, 0.6)$ at $t^*_1=0.3n$  which
                  changes to  $\theta^*_3=(1.0, 0.25)$\\ \hspace*{2.6cm} at $t^*_2=0.7n$ ($K^*=3$).
\end{itemize}

\begin{figure}[h!]
\begin{center}
\includegraphics[height=5.5cm, width=15.75cm]{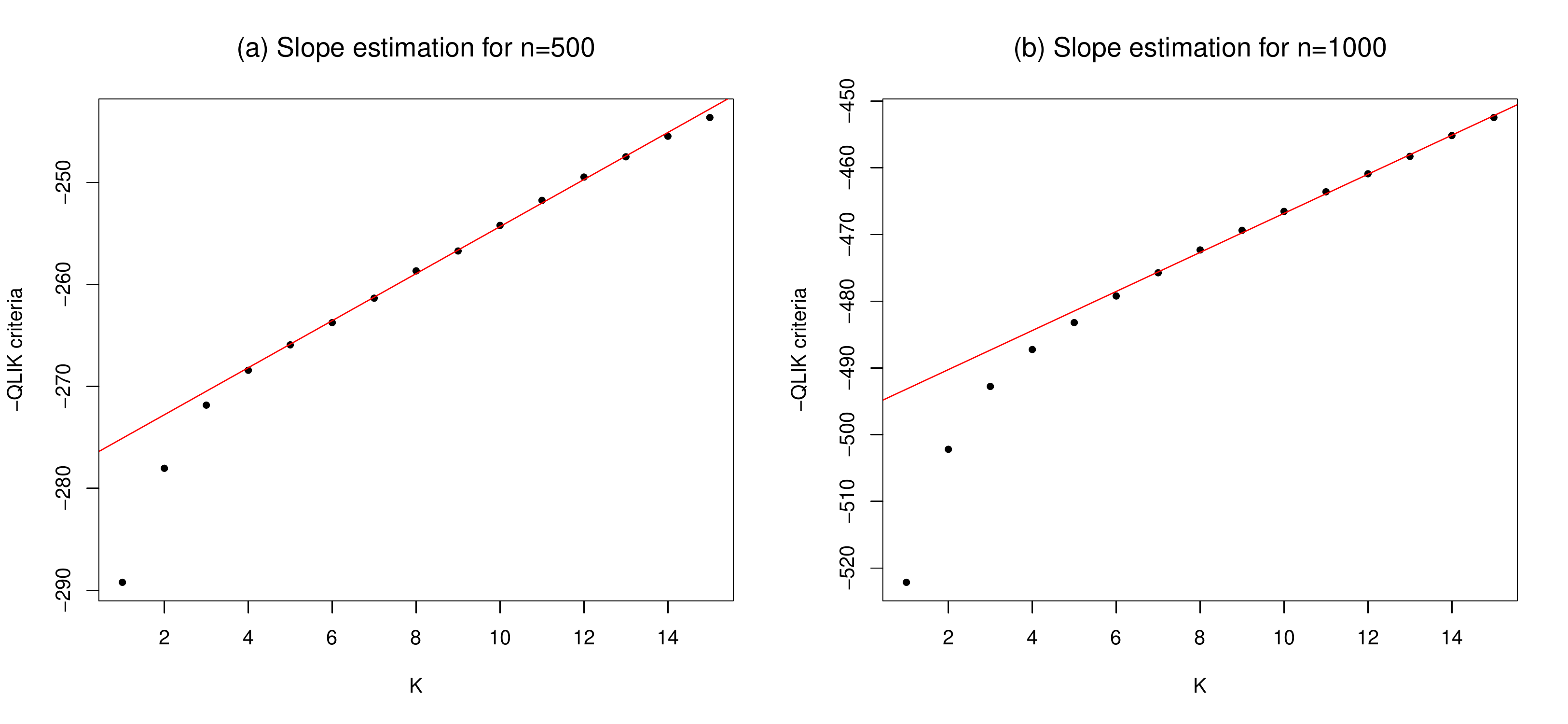}
\end{center}
\vspace{-1cm}
\caption{\it The curve of $-\min_{\underline{t},\underline{\theta}} QLIK(K)$, for $1 \leq K \leq K_{\max}$ for a Poisson-INARCH(1) process in scenario ${\bf IA}_2$.
The solid line represents the linear part of this curve with slope $\widehat{\kappa}_n/2 = 2.307$ when $n = 500$
and $\widehat{\kappa}_n/2  = 2.928$ when $n = 1000$.}
\label{Graphe_INARCH_1}
\end{figure}

\begin{figure}[h!]
\begin{center}
\includegraphics[height=5.5cm, width=15.75cm]{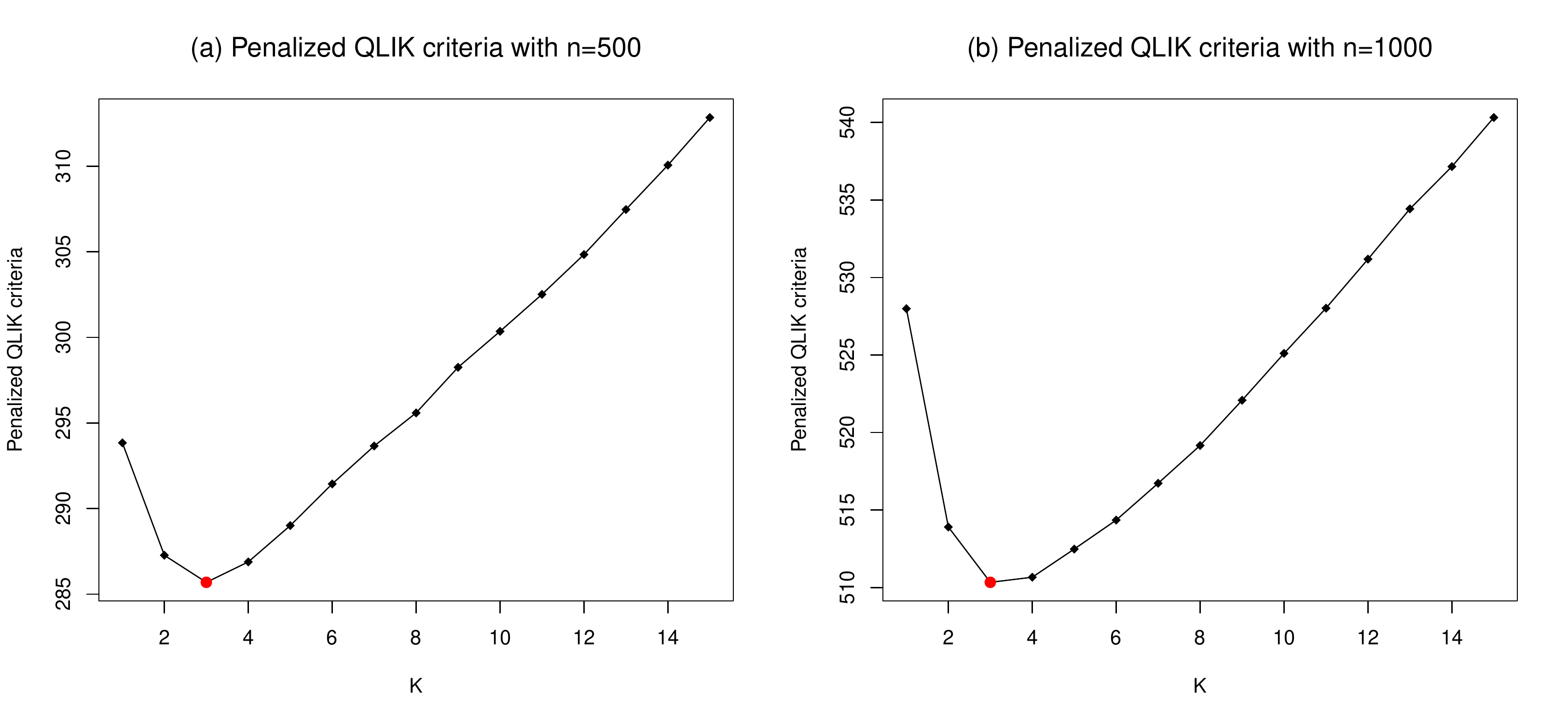}
\end{center}
\vspace{-1cm}
\caption{\it The graph $(K,\min_{\underline{t},\underline{\theta}} penQLIK(K))$, for $1 \leq K \leq K_{\max}$ for a Poisson-INARCH(1) process in scenario ${\bf IA}_2$.} \label{Graphe_INARCH_2}
\end{figure}

\begin{figure}[h!]
\begin{center}
\includegraphics[height=10.9cm, width=13cm]{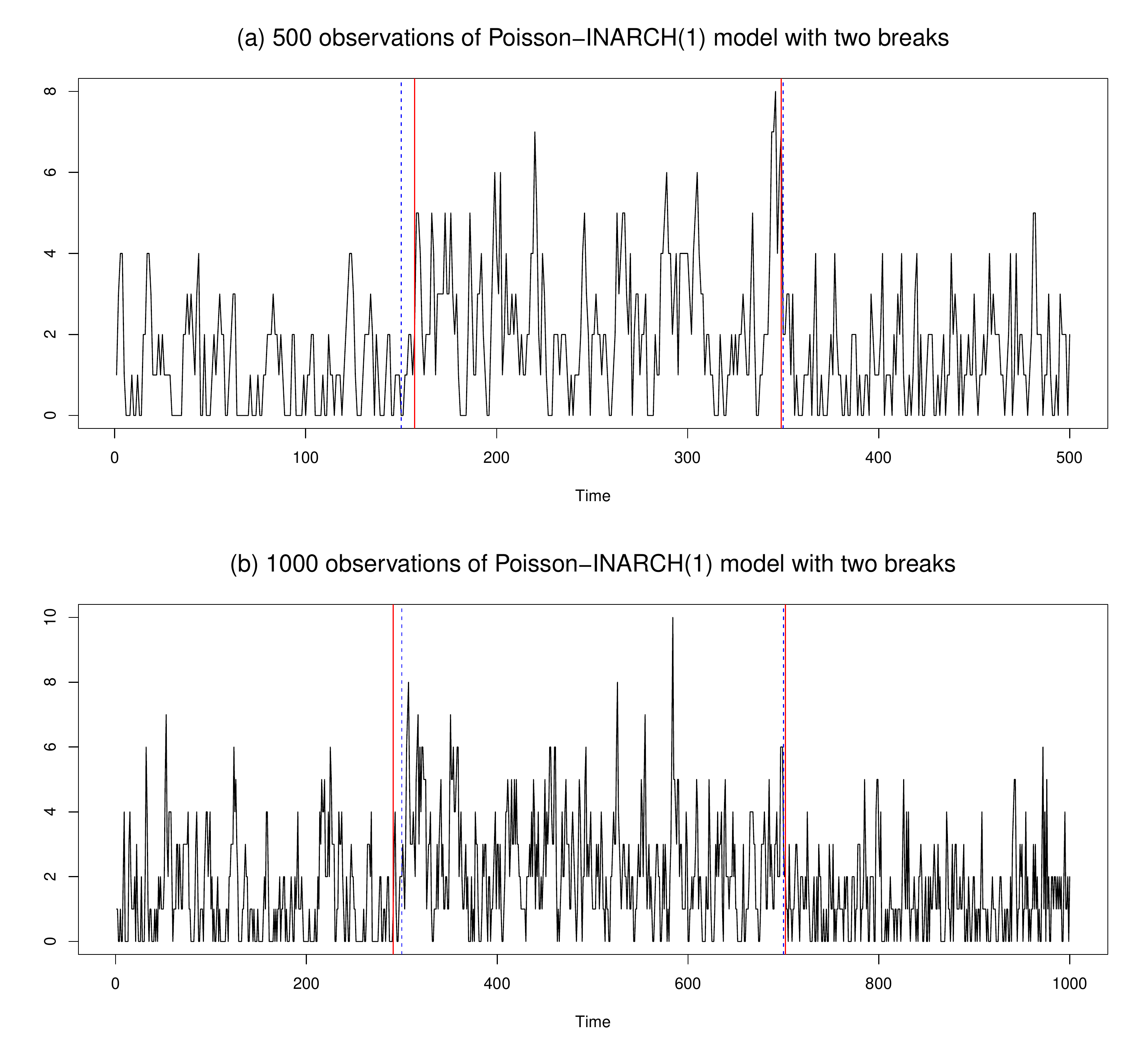}
\end{center}
\vspace{-.9cm}
\caption{\it The estimated of breakpoints for a trajectory of a Poisson-INARCH(1) process in scenario ${\bf IA}_2$. The
solid lines represent the estimated instants of breaks and the dotted lines represent the true ones.} \label{Graphe_INARCH_3}
\end{figure}

\noindent For scenario {\bf IA}$_2$,  Figure \ref{Graphe_INARCH_1} shows the slope of the linear part of the $-QLIK$ criteria minimized in $(\underline{t}, \underline{\theta})$. We obtain $\widehat{\kappa}_n \approx 4.6$ for $n=500$ and $\widehat{\kappa}_n \approx 5.9$ for $n=1000$.
Using these above values for $\kappa_n$ ({\it i.e.} $\kappa_n= \widehat{\kappa}_n$), we minimize the $penQLIK$ in $(K, \underline{t}, \underline{\theta})$, with $1 \leq K \leq K_{\max}$. Figure \ref{Graphe_INARCH_2} displays the points $(K,\min_{\underline{t},\underline{\theta}} penQLIK(K))$ for $1 \leq K \leq K_{\max}=15$.
 One can see that the estimated number of segments is $\widehat{K}_n=3$ for $n=500$ and $n=1000$.
 The estimated instants of breaks are $\widehat{\underline{t}}_n=(157,349)~(\underline{t}^*=(150,350))$ for $n=500$ and $\widehat{\underline{t}}_n=(291,702)~(\underline{t}^*=(300, 700))$ for $n=1000$ (see Figure \ref{Graphe_INARCH_3}).

 \medskip

 \medskip

\noindent Now, we are going to generate $100$ replications of a Poisson-INGARCH(1,1) process following the scenarios ${\bf IA}_0$-${\bf IA}_2$.
  Table \ref{Freq_Poisson_INARCH} indicates the frequencies of number of replications where $\widehat{K}_n =K^*$, $\widehat{K}_n <K^*$ and $\widehat{K}_n >K^*$,
   for the regularization parameter $\kappa_n= \widehat{\kappa}_n, ~ \log n, ~ n^{1/3}$.
 For the scenarios ${\bf IA}_1$ and ${\bf IA}_2$, we also consider the replications where the true number of breaks is achieved  (\textit{i.e.}    $\widehat{K}_n=K^*$) and we present some elementary statistics of the estimated instants of breaks (see Table \ref{Freq_Poisson_INARCH}).

\begin{table}[h!]
\scriptsize
\centering
\caption{\it Breaks estimated after 100 replications for Poisson-INARCH(1) process
following scenarios  ${\bf IA}_0$-${\bf IA}_2$. The first three columns show the frequencies of the estimation of the true, low and high number of breaks.
The last three columns give some elementary statistics of the  change-point locations when the true number of breaks is achieved.
}
\label{Freq_Poisson_INARCH}
\vspace{.2cm}
\begin{tabular}{cclccccccc}

\Xhline{.9pt}
& & &&&&& \\
& &&& \multicolumn{3} {c} {Frequencies} &  \multicolumn{2} {c}{Mean $\pm$ s.d.}   & {Mean}\\
Scenarios&&&&$\widehat{K}_n=K^*$&$\widehat{K}_n<K^*$&$\widehat{K}_n >K^*$ &$\widehat{\tau}_1$&$\widehat{\tau}_2$&$\left\|\underline{\widehat{\tau}}_n-\underline{\tau}^*\right\|$\\
\Xhline{.72pt}
\hline
{\bf IA}$_0$ & $n=500$& $\kappa_n=\widehat{\kappa}_n$         &&$0.72$&$0.00$&$0.28$&&&\\
 ($K^*=1$)            &        & $\kappa_n=\log n$            &&$0.76$&$0.00$&$0.24$&&&\\
                      &        & $\kappa_n=n^{1/3}$           &&$0.94$&$0.00$&$0.06$&&&\\

                      \cline{2-10}
                      & $n=1000$& $\kappa_n=\widehat{\kappa}_n$ &&$0.94$&$0.00$&$0.06$&&&\\
                      &        & $\kappa_n=\log n$              &&$0.90$&$0.00$&$0.10$&&& \\
                      &        & $\kappa_n=n^{1/3}$             &&$1.00$&$0.00$&$0.00$&&&\\

 \Xhline{0.9pt}

{\bf IA}$_1$ & $n=500$& $\kappa_n=\widehat{\kappa}_n$  &&$0.76$&$0.03$&$0.21$&$0.497 \pm 0.064$&&$0.038$\\
 ($K^*=2$)            &        & $\kappa_n=\log n$     &&$0.83$&$0.03$&$0.14$&$0.495 \pm 0.066$&&$0.040$\\
                      &        & $\kappa_n=n^{1/3}$    &&$0.87$&$0.09$&$0.04$&$0.495 \pm  0.064$&&$0.038$\\

           \cline{2-10}
                      & $n=1000$& $\kappa_n=\widehat{\kappa}_n$ &&$0.89$&$0.00$&$0.11$&$0.507 \pm 0.033$&&$0.019$\\
                      &        & $\kappa_n=\log n$              &&$0.87$&$0.00$&$0.13$&$0.507 \pm 0.034$&&$0.020$\\
                      &        & $\kappa_n=n^{1/3}$             &&$0.98$&$0.00$&$0.02$&$0.506 \pm 0.032$&&$0.019$\\

                       \Xhline{0.9pt}

{\bf IA}$_2$ &$n=500$ & $\kappa_n=\widehat{\kappa}_n$  &&$0.62$&$0.13$&$0.25$&$0.311 \pm  0.071$&$0.689 \pm 0.060$&$ 0.061$\\
 ($K^*=3$)            &        & $\kappa_n=\log n$     &&$0.73$&$0.12$&$0.15$&$0.317 \pm  0.073$&$0.690 \pm 0.072 $&$ 0.067$\\
                      &        & $\kappa_n=n^{1/3}$    &&$0.64$&$0.33$&$0.03$&$0.310 \pm  0.058$&$0.685 \pm 0.070 $&$ 0.061$\\
          \cline{2-10}
               &$n=1000$& $\kappa_n=\widehat{\kappa}_n$  &&$0.87$&$0.00$&$0.13$&$0.300 \pm 0.034$&$0.693 \pm 0.030$&$0.034$\\
               &        & $\kappa_n=\log n$              &&$0.84$&$0.00$&$0.16$&$0.302 \pm 0.043$&$  0.692\pm 0.030$&$0.038$\\
               &        & $\kappa_n=n^{1/3}$             &&$0.93$&$0.05$&$0.02$&$0.300\pm  0.051$&$0.694 \pm 0.028 $&$0.037$\\

 \Xhline{.9pt}
\end{tabular}
\end{table} ~\\

\noindent The results of Table \ref{Freq_Poisson_INARCH} show that for the penalties considered, the performance increase with $n$ in all scenarios.
 In accordance with Theorem \ref{th1}, the consistency of the penalties $\log n$ and $n^{1/3}$ is numerically convincing.
 Moreover, the $n^{1/3}$ penalty outperforms the other procedures when $n=1000$.

\subsubsection{Poisson-INGARCH models}
We consider the problem (\ref{model_pwc}) for a Poisson-INGARCH(1,1), {\it i.e.}  $(Y_1,\cdots,Y_n)$ is a trajectory of the process  $Y=\{Y_{t},t\in \Z \}$ satisfying:
\begin{equation}\label{Poisson_pwc_INGARCH}
    Y_{t}|\mathcal{F}_{t-1} \sim \mathcal{P} (\lambda_{t})~;~  \lambda_{t}= 
    \alpha^{(j)}_0 + \alpha^{(j)} Y_{t-1}+ \beta^{(j)} \lambda_{t-1}, ~  ~  \forall~ t \in T^*_j,~~   \forall~ j \in \left\{1,\cdots,K^*\right\}.
   \end{equation}
   The parameter vector  is $\theta^*_j=(\alpha^{(j)}_0 ,\alpha^{(j)}, \beta^{(j)})$, for all $j \in \{1,\cdots,K^*\}$.

  \medskip

   \medskip

\noindent For $n=500$ and $n=1000$, we generate $100$ replications of the model (\ref{Poisson_pwc_INGARCH})  in the following situations:
\begin{itemize}
    \item {\bf scenario IG}$_0$: $\theta^{*}_1=(1.0,0.2,0.15)$ is constant ($K^*=1$) ;
    \item {\bf scenario IG}$_1$: $\theta^*_1=(1.0,0.2,0.15)$ changes to $\theta^*_2=(1.0,0.45,0.15)$ at $t^*=0.5n$ ($K^*=2$) ;
    \item {\bf scenario IG}$_2$: $\theta^*_1=(0.1,0.3,0.6)$ changes to $\theta^*_2=(0.5,0.3,0.6)$ at $t^*_1=0.3n$ which
                  changes to\\ \hspace*{2.6cm} $\theta^*_3=(0.5,0.3,0.2)$ at $t^*_2=0.7n$ ($K^*=3$).
\end{itemize}

\begin{table}[h!]
\scriptsize
\centering
\caption{\it Breaks estimated after 100 replications for Poisson-INGARCH(1,1) process
following scenarios  ${\bf IG}_0$-${\bf IG}_2$. The first three columns show the frequencies of the estimation of the true, low and high number of breaks.
The last three columns give some elementary statistics of the change-point locations when the true number of breaks is achieved.}
\label{Freq_Poisson_INGARCH}
\vspace{.2cm}
\begin{tabular}{cclccccccc}

\Xhline{.9pt}
& & &&&&& \\
& &&& \multicolumn{3} {c} {Frequencies} &  \multicolumn{2} {c}{Mean $\pm$ s.d.}   & {Mean}\\
Scenarios&&&&$\widehat{K}_n=K^*$&$\widehat{K}_n<K^*$&$\widehat{K}_n >K^*$ &$\widehat{\tau}_1$&$\widehat{\tau}_2$&$\left\|\underline{\widehat{\tau}}_n-\underline{\tau}^*\right\|$\\
\Xhline{.72pt}
\hline
{\bf  IG}$_0$ & $n=500$& $\kappa_n=\widehat{\kappa}_n$          &&$0.86$&$0.00$&$0.14$&&&\\
 ($K^*=1$)            &        & $\kappa_n=\log n$              &&$0.96$&$0.00$&$0.04$&&&\\
                      &        & $\kappa_n=n^{1/3}$             &&$1.00$&$0.00$&$0.00$&&&\\

      \cline{2-10}
                      & $n=1000$& $\kappa_n=\widehat{\kappa}_n$ &&$0.92$&$0.00$&$0.08$&&&\\
                      &        & $\kappa_n=\log n$              &&$0.96$&$0.00$&$0.04$&&&\\
                      &        & $\kappa_n=n^{1/3}$             &&$1.00$&$0.00$&$0.00$&&&\\

 \Xhline{.9pt}
{\bf IG}$_1$ & $n=500$& $\kappa_n=\widehat{\kappa}_n$           &&$0.75$&$0.05$&$0.20$&$0.515 \pm 0.066$&&$0.038$\\
 ($K^*=2$)            &        & $\kappa_n=\log n$              &&$0.70$&$0.03$&$0.27$&$0.514 \pm 0.073$&&$0.040$\\
                      &        & $\kappa_n=n^{1/3}$             &&$0.78$&$0.06$&$0.16$&$0.512 \pm 0.066$&&$0.038$\\

           \cline{2-10}
                      & $n=1000$& $\kappa_n=\widehat{\kappa}_n$ &&$0.75$&$0.05$&$0.20$&$0.507\pm 0.031$&&$0.019$\\
                      &        & $\kappa_n=\log n$              &&$0.58$&$0.00$&$0.42$&$0.508\pm 0.034$&&$0.021$\\
                      &        & $\kappa_n=n^{1/3}$             &&$0.83$&$0.03$&$0.13$&$0.501\pm 0.048$&&$0.022$\\

                       \Xhline{.9pt}

{\bf  IG}$_2$ &$n=500$ & $\kappa_n=\widehat{\kappa}_n$        &&$0.53$&$0.41$&$0.06$&$0.299\pm 0.078$&$0.691 \pm 0.073$&$ 0.053$\\
 ($K^*=3$)            &        & $\kappa_n=\log n$            &&$0.58$&$0.23$&$0.19$&$0.299 \pm 0.074$&$ 0.693\pm 0.070$&$0.049$\\
                      &        & $\kappa_n=n^{1/3}$           &&$0.37$&$0.49$&$0.14$&$0.300 \pm 0.076$&$0.697 \pm 0.015$&$0.047$\\

          \cline{2-10}
                    &$n=1000$& $\kappa_n=\widehat{\kappa}_n$  &&$0.62$&$0.26$&$0.12$&$0.293 \pm 0.050$&$0.702 \pm 0.010$&$0.025$\\
                    &        & $\kappa_n=\log n$              &&$0.60$&$0.06$&$0.34$&$ 0.293\pm 0.051$&$0.702 \pm 0.010$&$0.026$\\
                    &        & $\kappa_n=n^{1/3}$             &&$0.56$&$0.29$&$0.15$&$0.301\pm  0.029$&$0.699 \pm 0.011$&$0.016$\\

 \Xhline{.9pt}
\end{tabular}
\end{table}
\noindent Table \ref{Freq_Poisson_INGARCH} indicates the frequencies of the true number of breaks estimated and some elementary statistics of the estimators of the change-point locations.
 It appears that the results of the  $n^{1/3}$-penalty and the slope procedure are quite satisfactory except for the case of two breaks.
In this later case, the $n^{1/3}$-penalty and the slope procedure over-penalizes the number of breaks, while the $\log n$-penalty under-penalizes.
But, overall, the performances of the proposed procedures increase with $n$ and the estimation of the breakpoints locations is well achieved.


\subsubsection{Negative binomial INGARCH models}
We consider the problem (\ref{model_pwc}) for a  negative binomial INGARCH(1,1) (NB-INGARCH(1,1)), {\it i.e.}  $(Y_1,\cdots,Y_n)$ is a trajectory of the process  $Y=\{Y_{t},t\in \Z \}$ satisfying:
\begin{equation}\label{NB_pwc_INGARCH}
    Y_{t}|\mathcal{F}_{t-1} \sim NB(r,p_{t})~;~ r\frac{(1-p_{t})}{p_{t}}= \lambda_{t}= 
    \alpha^{(j)}_0 + \alpha^{(j)} Y_{t-1}+ \beta^{(j)} \lambda_{t-1}, ~  ~  \forall~ t \in T^*_j,~~   \forall~ j \in \left\{1,\cdots,K^*\right\};
   \end{equation}
  where the parameter vector  is $\theta^*_j=(\alpha^{(j)}_0 ,\alpha^{(j)}, \beta^{(j)})$, for all $j \in \{1,\cdots,K^*\}$ and $NB(r, p)$ denotes the negative binomial distribution with parameters $r$ and $p$.\\
\noindent For $r=14$ (used for  transaction data, see \cite{Diop2017}), $n=500$ and $n=1000$, we generate a sample $(Y_1, \cdots ,Y_n)$ in the following situations:
\begin{itemize}
    \item {\bf scenario NB-IG}$_0$: $\theta^{*}_1=(1.0,0.2,0.15)$ is constant ($K^*=1$) ;
    \item {\bf scenario NB-IG}$_1$: $\theta^*_1=(1,0.2,0.15)$ changes to $\theta^*_2=(1,0.45,0.15)$ at $t^*=0.5n$ ($K^*=2$) ;
    \item {\bf scenario NB-IG}$_2$: $\theta^*_1=(0.1,0.3,0.6)$ changes to $\theta^*_2=(0.5,0.3,0.6)$ at $t^*_1=0.3n$ which
                  changes to \\
                   \hspace*{3.2cm} $\theta^*_3=(0.5,0.3,0.2)$  at $t^*_2=0.7n$ ($K^*=3$).
\end{itemize}
  \begin{table}[h!]
\scriptsize
\centering
\caption{\it Breaks estimated after 100 replications for NB-INGARCH(1,1) process following scenarios  ${\bf NB-IG}_0$-${\bf NB-IG}_2$.
 The first three columns show the frequencies of the estimation of the true, low and high number of breaks.
The last three columns give some elementary statistics of the  change-point locations when the true number of breaks is achieved.}
\label{Freq_NB_INGARCH}
\vspace{.2cm}
\begin{tabular}{cclccccccc}

\Xhline{.9pt}
& & &&&&& \\
& &&& \multicolumn{3} {c} {Frequencies} &  \multicolumn{2} {c}{Mean $\pm$ s.d.}   & {Mean}\\
Scenarios &&&&$\widehat{K}_n=K^*$&$\widehat{K}_n<K^*$&$\widehat{K}_n >K^*$ &$\widehat{\tau}_1$&$\widehat{\tau}_2$&$\left\|\underline{\widehat{\tau}}_n-\underline{\tau}^*\right\|$\\
\Xhline{.72pt}
\hline
{\bf NB-IG}$_0$ & $n=500$& $\kappa_n=\widehat{\kappa}_n$  &&$0.90$&$0.00$&$0.10$&&&\\
 ($K^*=1$)            &        & $\kappa_n=\log n$                 &&$0.95$&$0.00$&$0.05$&&&\\
                      &        & $\kappa_n=n^{1/3}$                &&$0.98$&$0.00$&$0.02$&&&\\
                      \cline{2-10}
                      & $n=1000$& $\kappa_n=\widehat{\kappa}_n$ &&$0.92$&$0.00$&$0.08$&&&\\
                      &        & $\kappa_n=\log n$              &&$0.94$&$0.00$&$0.06$&&&\\
                      &        & $\kappa_n=n^{1/3}$             &&$0.98$&$0.00$&$0.02$&&&\\

 \Xhline{.9pt}
{\bf NB-IG}$_1$ & $n=500$& $\kappa_n=\widehat{\kappa}_n$   &&$0.60$&$0.09$&$0.31$&$0.512 \pm 0.122$&&$0.072$\\
 ($K^*=2$)            &        & $\kappa_n=\log n$         &&$0.55$&$0.04$&$0.41$&$0.514 \pm 0.110$&&$0.063$\\
                      &        & $\kappa_n=n^{1/3}$        &&$0.65$&$0.12$&$0.23$&$0.519 \pm  0.106$&&$0.060$\\

           \cline{2-10}
                      & $n=1000$& $\kappa_n=\widehat{\kappa}_n$ &&$0.69$&$0.02$&$0.29$&$ 0.507\pm 0.065$&& $0.037$\\
                      &        & $\kappa_n=\log n$              &&$0.56$&$0.00$&$0.44$&$0.500\pm 0.054$&&  $0.030$\\
                      &        & $\kappa_n=n^{1/3}$             &&$0.83$&$0.02$&$0.15$&$0.505 \pm 0.061$&& $0.037$\\

                       \Xhline{.9pt}

{\bf NB-IG}$_2$ &$n=500$ & $\kappa_n=\widehat{\kappa}_n$  &&$0.45$&$0.51$&$0.04$&$0.330\pm 0.084$&$ 0.696 \pm 0.040$&$ 0.057$\\
 ($K^*=3$)            &        & $\kappa_n=\log n$        &&$0.41$&$0.13$&$0.46$&$0.319\pm 0.080$&$ 0.700 \pm 0.026$&$0.057$\\
                      &        & $\kappa_n=n^{1/3}$       &&$0.43$&$0.28$&$0.29$&$0.328\pm 0.061$&$0.685 \pm 0.066 $&$0.055$\\

          \cline{2-10}
                      &$n=1000$& $\kappa_n=\widehat{\kappa}_n$  &&$0.70$&$0.25$&$0.05$&$0.304\pm 0.060$&$0.699 \pm0.020$&$0.033$\\
                      &        & $\kappa_n=\log n$              &&$0.54$&$0.02$&$0.44$&$0.299\pm 0.066$&$0.699\pm 0.014$&$0.033$\\
                      &        & $\kappa_n=n^{1/3}$             &&$0.68$&$0.19$&$0.13$&$0.325\pm 0.090$&$0.697\pm 0.018$&$0.043$\\

 \Xhline{.9pt}
\end{tabular}
\end{table}
Once again, it appears in Table \ref{Freq_NB_INGARCH} that the performances of the proposed procedures increase with $n$
and the estimation of the breakpoints locations remain satisfactory even in this case where the Poisson quasi-likelihood used is quite different from the true distribution
of the observations.


\subsubsection{Binary Time Series}

Consider the problem (\ref{model_pwc}) for a binary INARCH(1) (BIN-INARCH(1)) time series model, {\it i.e.}  $(Y_1,\cdots,Y_n)$ is a trajectory of the process  $Y=\{Y_{t},t\in \Z \}$ satisfying:
\begin{equation}\label{BIN_pwc_INGARCH}
    Y_{t}|\mathcal{F}_{t-1} \sim \mathcal{B}(p_{t})~;~ p_{t} = \lambda_{t}=
    \alpha^{(j)}_0 + \alpha^{(j)} Y_{t-1}, ~  ~  \forall~ t \in T^*_j,~~   \forall~ j \in \left\{1,\cdots,K^*\right\};
   \end{equation}
  where the parameter vector  is $\theta^*_j=(\alpha^{(j)}_0 ,\alpha^{(j)})$ (with $0<\alpha^{(j)}_0 +\alpha^{(j)}<1$), for all $j \in \{1,\cdots,K^*\}$ and $\mathcal{B}(p)$ denotes the Bernoulli distribution with parameter $p$.\\
 \noindent For $n=500$ and $n=1000$, we generate a sample $(Y_1, \cdots ,Y_n)$ in the following situations:
\begin{itemize}
    \item {\bf scenario BIN-IA}$_0$: $\theta^{*}_1=(0.15,0.75)$ is constant ($K^*=1$) ;
    \item {\bf scenario BIN-IA}$_1$: $\theta^*_1=(0.15,0.75)$ changes to $\theta^*_2=(0.04,0.60)$ at $t^*=0.5n$ ($K^*=2$) ;
    \item {\bf scenario BIN-IA}$_2$: $\theta^*_1=(0.15,0.75)$ changes to $\theta^*_2=(0.04,0.60)$ at $t^*_1=0.3n$ which
                  changes to \\ \hspace*{3.4cm}  $\theta^*_3=(0.25,0.35)$ at $t^*_2=0.7n$ ($K^*=3$).
\end{itemize}
 \noindent The scenario \textbf{BIN-IA}$_1$ is related and close to the real data example (see below).

\begin{table}[h!]
\scriptsize
\centering
\caption{ \it Breaks estimated after 100 replications for BIN-INARCH(1) process following scenarios  ${\bf BIN-IA}_0$-${\bf BIN-IA}_2$.
 The first three columns show the frequencies of the estimation of the true, low and high number of breaks.
The last three columns give some elementary statistics of the  change-point locations when the true number of breaks is achieved.}
\label{Freq_BIN_INARCH}
\vspace{.2cm}
\begin{tabular}{cclccccccc}

\Xhline{.9pt}
& & &&&&& \\
& &&& \multicolumn{3} {c} {Frequencies} &  \multicolumn{2} {c}{Mean $\pm$ s.d.}   & {Mean}\\
Scenarios&&&&$\widehat{K}_n=K^*$&$\widehat{K}_n<K^*$&$\widehat{K}_n >K^*$ &$\widehat{\tau}_1$&$\widehat{\tau}_2$&$\left\|\underline{\widehat{\tau}}_n-\underline{\tau}^*\right\|$\\
\Xhline{.72pt}
\hline
{\bf BIN-IA}$_0$ & $n=500$& $\kappa_n=\widehat{\kappa}_n$  &&$0.84$&$0.00$&$0.16$&&&\\
 ($K^*=1$)            &        & $\kappa_n=\log n$         &&$0.98$&$0.00$&$0.02$&&&\\
                      &        & $\kappa_n=n^{1/3}$        &&$0.98$&$0.00$&$0.02$&&&\\

                      \cline{2-10}
                      & $n=1000$& $\kappa_n=\widehat{\kappa}_n$ &&$0.86$&$0.00$&$0.14$&&&\\
                      &        & $\kappa_n=\log n$              &&$1.00$&$0.00$&$0.00$&&&\\
                      &        & $\kappa_n=n^{1/3}$             &&$1.00$&$0.00$&$0.00$&&&\\

 \Xhline{.9pt}
{\bf BIN-IA}$_1$ & $n=500$& $\kappa_n=\widehat{\kappa}_n$  &&$0.69$&$0.10$&$0.21$&$0.499 \pm  0.091$&&$0.055$\\
 ($K^*=2$)            &        & $\kappa_n=\log n$         &&$0.72$&$0.27$&$0.01$&$ 0.491 \pm  0.087$&&$0.051$\\
                      &        & $\kappa_n=n^{1/3}$        &&$0.52$&$0.48$&$0.00$&$0.484 \pm 0.091  $&&$0.054$\\

           \cline{2-10}
                      & $n=1000$& $\kappa_n=\widehat{\kappa}_n$ &&$0.89$&$0.00$&$0.11$&$0.499 \pm  0.036$&&$0.020$\\
                      &        & $\kappa_n=\log n$              &&$0.96$&$0.01$&$0.03$&$0.500 \pm  0.035$&&$0.019$\\
                      &        & $\kappa_n=n^{1/3}$             &&$0.85$&$0.15$&$0.00$&$0.484 \pm 0.091  $&&$0.054$\\

                       \Xhline{.9pt}

{\bf BIN-IA}$_2$ &$n=500$ & $\kappa_n=\widehat{\kappa}_n$ &&$0.75$&$0.18$&$0.07$&$0.324\pm  0.094$&$0.695 \pm 0.044 $&$ 0.060$\\
 ($K^*=3$)            &        & $\kappa_n=\log n$        &&$0.43$&$0.57$&$0.00$&$0.312\pm  0.044$&$0.694 \pm 0.030 $&$ 0.035$\\
                      &        & $\kappa_n=n^{1/3}$       &&$0.23$&$0.77$&$0.00$&$0.306\pm  0.033$&$0.702 \pm 0.017 $&$ 0.026$\\

          \cline{2-10}
                      &$n=1000$& $\kappa_n=\widehat{\kappa}_n$  &&$0.95$&$0.03$&$0.02$&$0.303\pm 0.046$&$0.696\pm 0.019$&$0.028$\\
                      &        & $\kappa_n=\log n$              &&$0.90$&$0.10$&$0.00$&$0.299\pm 0.037$&$0.697\pm 0.017$&$0.025$\\
                      &        & $\kappa_n=n^{1/3}$             &&$0.52$&$0.48$&$0.00$&$0.296\pm 0.021$&$0.697\pm 0.019$&$0.020$\\

 \Xhline{.9pt}
\end{tabular}
\end{table}

\noindent Table \ref{Freq_BIN_INARCH} shows that the procedure provides satisfactory results with BIN-INARCH(1) model, except that the $n^{1/3}$-penalty in the case of two breaks.
But, the performances of these procedures increase with $n$ and the breakpoints locations are overall well estimated.


\subsubsection{INARCH$(\infty)$ models}
Now, consider a Poisson-INARCH$(\infty)$, {\it i.e.}  $(Y_1,\cdots,Y_n)$ is a trajectory of the process  $Y=\{Y_{t},t\in \Z \}$ satisfying:
\begin{equation}\label{Poisson_pwc_INARCH_Inf}
    Y_{t}|\mathcal{F}_{t-1} \sim \mathcal{P}(\lambda_{t})~~;~~  \lambda_{t}=
    \alpha^{(j)}_{0} + \sum_{k =1}^{\infty}\alpha_{k} Y_{t-k} , ~~  ~  \forall~ t \in T^*_j,~~   \forall~ j \in \left\{1,\cdots,K^*\right\},
   \end{equation}
 where $\alpha^{(j)}_{0}>0$, $\alpha_{k} \geq 0$ (for all $k \geq 1$ and $j =1,\cdots,K^*$) and $\underset{k \geq 1}{\sum}\alpha_{k}<1$;
 that is, we focus on the change in the parameter $\alpha^{(j)}_{0}$.
 This process corresponds to a particular case of the problem (\ref{model_pwc}) with $f(y_1,y_2,\cdots,\alpha^{(j)}_{0})= \alpha^{(j)}_0 + \sum_{k =1}^{\infty}\alpha_k y_{k}$
  for each regime $j \in \{1,\cdots,K^*\}$.
 We deal with a scenario where the consistency of the BIC procedure is not ensured.
 Therefore, we consider the Riemanian case with $ \alpha_{k} =O(k^{-1.7}) $ (in the scenario detailed below).
  More precisely, we consider the model (\ref{Poisson_pwc_INARCH_Inf}) with
  \[  \lambda_{t}= \alpha^{(j)}_{0} + \frac{1}{2.2} \sum_{k =1}^{\infty} \frac{1}{k^{1.7}} Y_{t-k}  .\]
  The number $1/2.2$ is obtained from the values of the Riemann zeta function, and allows the condition $\frac{1}{2.2} \sum_{k =1}^{\infty} \frac{1}{k^{1.7}} < 1$.
 According to Theorem \ref{th1} and Remark \ref{rmk1}, if the regularization parameter verifies
 $\kappa_n = O(n^{\delta})$ with $\delta > 0.3$, then the consistency holds. Thus, the consistency of the BIC penalty is not ensured.


\medskip

\noindent Now, for $n=500$ and $n=100$, we generate a trajectory $(Y_1,\cdots,Y_n)$ of the model (\ref{Poisson_pwc_INARCH_Inf}) in the following scenarios:
\begin{itemize}
    \item {\bf scenario IA-INF}$_0$: $\alpha^{(1)}_{0}=0.5$ is constant ($K^*=1$) ;
    \item {\bf scenario IA-INF}$_1$: $\alpha^{(1)}_{0}=0.5$ changes to $\alpha^{(2)}_{0}=0.1$ at $t^*=0.5n$ ($K^*=2$).
\end{itemize}
\begin{table}[h!]
\scriptsize
\centering
\caption{\it Breaks estimated after 100 replications for INARCH$(\infty)$ process following scenarios  ${\bf IA-INF}_0$ and ${\bf IA-INF}_1$.
 The first three columns show the frequencies of the estimation of the true, low and high number of breaks.
The last two columns give some elementary statistics of the  change-point locations when the true number of breaks is achieved.}
\label{Freq_IA-INF}
\vspace{.2cm}
\begin{tabular}{cclcccccc}

\Xhline{.9pt}
& & &&&&& \\
& &&& \multicolumn{3} {c} {Frequencies} & {Mean $\pm$ s.d.}  & {Mean}\\
Scenarios&&&&$\widehat{K}_n=K^*$&$\widehat{K}_n<K^*$&$\widehat{K}_n >K^*$ &$\widehat{\tau}_1$&$\left\|\underline{\widehat{\tau}}_n-\underline{\tau}^*\right\|$\\
\Xhline{.72pt}
\hline
{\bf IA-INF}$_0$ & $n=500$& $\kappa_n=\widehat{\kappa}_n$  &&$0.15$&$0.00$&$0.85$&&\\
 ($K^*=1$)            &        & $\kappa_n=\log n$         &&$0.56$&$0.00$&$0.44$&&\\
                      &        & $\kappa_n=n^{1/3}$        &&$0.84$&$0.00$&$0.16$&&\\

                      \cline{2-9}
                      & $n=1000$& $\kappa_n=\widehat{\kappa}_n$ &&$0.17$&$0.00$&$0.83$&&\\
                      &        & $\kappa_n=\log n$              &&$0.57$&$0.00$&$0.43$&&\\
                      &        & $\kappa_n=n^{1/3}$             &&$0.95$&$0.00$&$0.05$&&\\

 \Xhline{.9pt}
{\bf IA-INF}$_1$ & $n=500$& $\kappa_n=\widehat{\kappa}_n$  &&$0.66$&$0.04$&$0.34$&$0.498 \pm   0.0057$&$0.003$\\\
 ($K^*=2$)            &        & $\kappa_n=\log n$         &&$0.49$&$0.02$&$0.49$&$0.498\pm  0.0013$&$0.002$\\
                      &        & $\kappa_n=n^{1/3}$        &&$0.82$&$0.03$&$0.15$&$0.497 \pm  0.0062$&$0.003$\\

           \cline{2-9}
                      & $n=1000$& $\kappa_n=\widehat{\kappa}_n$ &&$0.79$&$0.00$&$0.21$&$ 0.499\pm  0.0004$&$0.001$\\
                      &        & $\kappa_n=\log n$              &&$0.43$&$0.00$&$0.57$&$0.499 \pm  0.0006$&$0.001$\\
                      &        & $\kappa_n=n^{1/3}$             &&$0.93$&$0.00$&$0.07$&$ 0.499 \pm  0.0005$&$0.001$\\

 \Xhline{.9pt}
\end{tabular}
\end{table}
\noindent
In Table \ref{Freq_IA-INF}, on can see that the $n^{1/3}$-penalty uniformly  outperforms the other two procedures.
Moreover, the performances of the proposed procedures increase with $n$, except the $\log n$-penalty whose the performances decrease with $n$.
Hence, the consistency of the BIC procedure is quite questionable in this case.
%



\section{Real data application}

We apply our change-point procedure to two examples of real data series. To compute the estimator $\widehat{K}_n$, the $\widehat{\kappa}_n$
penalty is used with $u_n=\left[(\log (n))^{\delta} \right]$ (where $3/2\leq \delta \leq 2$) and $K_{\max}=15$.

\subsection{The US recession data}
\begin{figure}[h!]
\begin{center}
\includegraphics[height=5.5cm, width=11.5cm]{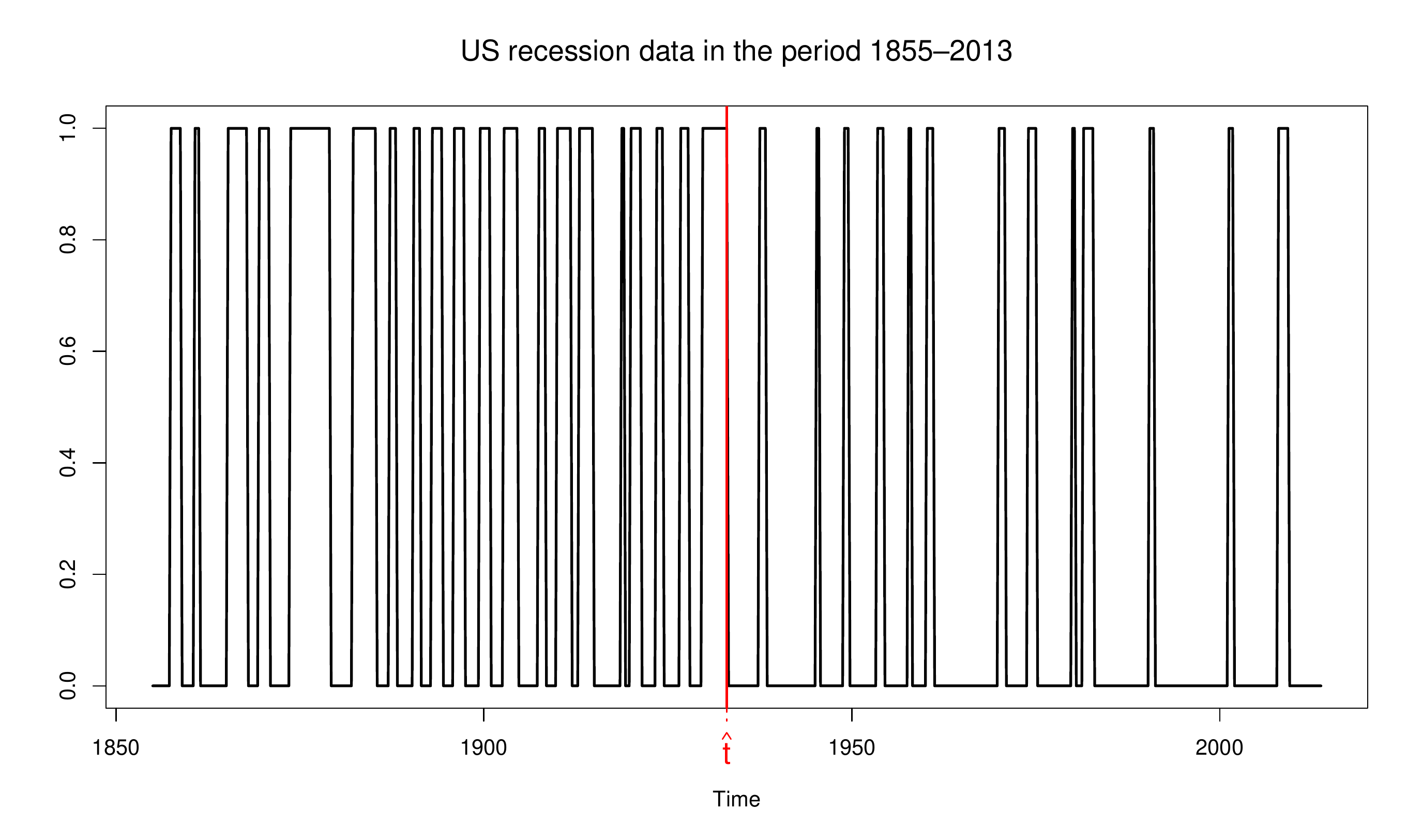}
\end{center}
\vspace{-1cm}
\caption{ The US recession data with the estimated location of the breakpoint $\widehat{t}$.}
\label{Result_US_recession_2}
\end{figure}

Firstly, we consider the series of the quarterly recession data from the USA for the period 1855-2013 (see Figure \ref{Result_US_recession_2}).
 This series $(Y_t)$ represents a binary variable that is equal to $1$ if there is a recession in at least $1$ month in the quarter and $0$ otherwise.
 There are 636 quarterly observations obtained from The National Bureau of Economic Research.
These data have already been analyzed by  several authors.
 Hudecov{\'a} (2013) has applied a change-point procedure based on a normalized cumulative sums of residuals and
has found a break in the first quarter of 1933. Recently, Diop and Kengne (2017) have applied a change-point test based on the maximum likelihood estimator of the model's
parameter and have detected a break in the last quarter of 1932.
\medskip

\begin{figure}[h!]
\begin{center}
\includegraphics[height=5.5cm, width=15cm]{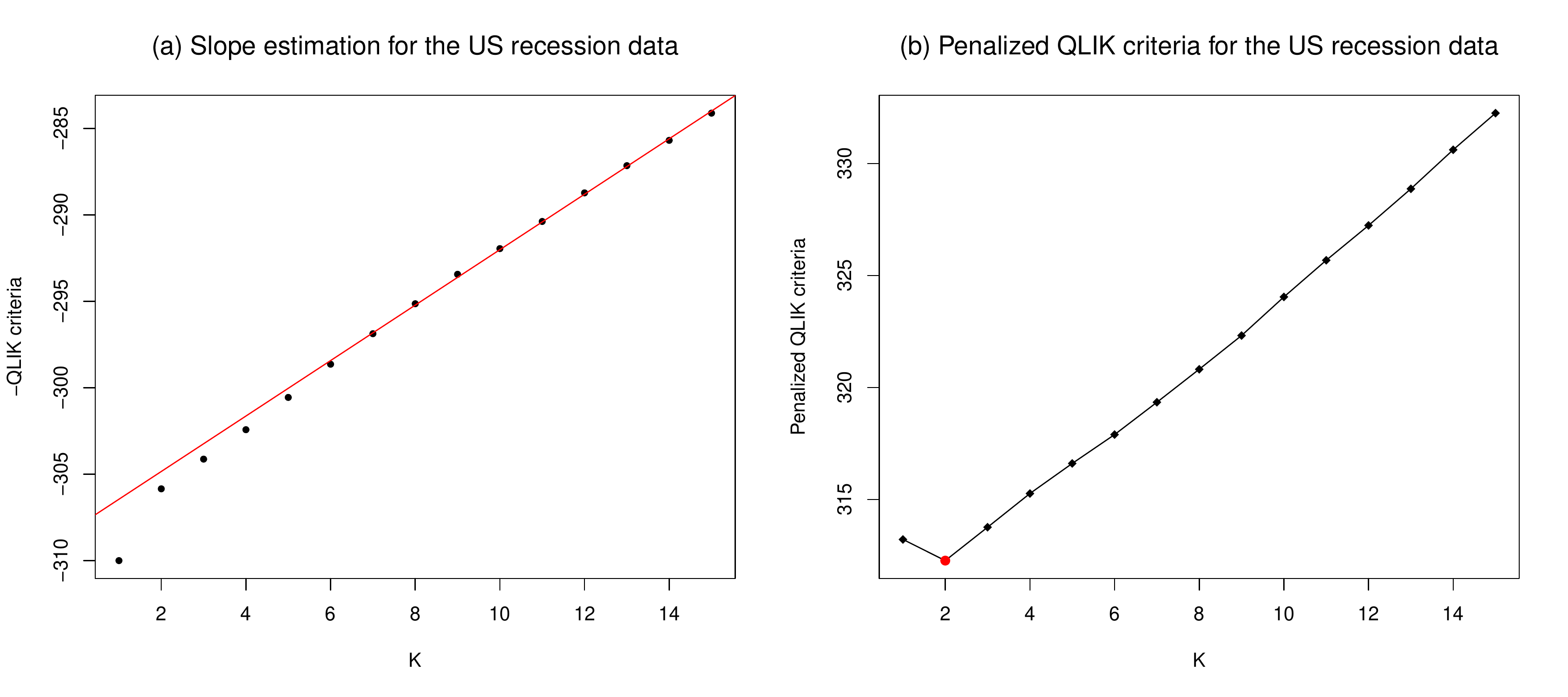}
\end{center}
\vspace{-.8cm}
\caption{ The curve of $-\min_{\underline{t},\underline{\theta}} QLIK(K)$ and the graph $(K,\min_{\underline{t},\underline{\theta}} penQLIK(K))$ for the US recession data with a INARCH(1) model.
The solid line represents the linear part of this curve with slope $\widehat{\kappa}_n/2 = 1.605$.}
\label{Result_US_recession_1}
\end{figure}

\noindent We consider the INARCH(1) representation and apply the $penQLIK$ contrast procedure. This choice is motivated by the fact that the estimation of the last component
of $\theta$ ({\it i.e.} the parameter $\beta$) is not significant in the INGARCH(1,1) representation (see Diop and Kengne \cite{Diop2017}).
The test of nullity of one coefficient (TNOC) proposed by Ahmad and Francq \cite{Francq2016}, applied a posteriori (after change-point detection)  
also confirms these results. 
 As noted in the implementation of the dynamic programming algorithm, we begin by the calibration of the regularization parameter $\kappa_n$.
 The slope estimation procedure applied with $u_n=\left[(\log n)^{2} \right]$ returns the values $\widehat{\kappa}_n \approx 3.21$
 and the estimation of the number of segments is $\widehat{K}_n=2$, {\it i.e.} one break is detected (see Figure \ref{Result_US_recession_1}).
 The location of the breakpoint estimated  is $\widehat{t}=313$. The change detected at $\widehat{t}=313$ corresponds to the first quarter of 1933
 (see Figure \ref{Result_US_recession_2}). These results are in concordance with those obtained by Diop and Kengne (2017)  and  Hudecov{\'a} (2013).
  The estimated model with one breakpoint is
\[\E(Y_t | \mathcal{F}_{t-1})=\lambda_{t}=\left\{
\begin{array}{l}
  0.120 + 0.749 Y_{t-1} , \text{ for } t \leq 313 \\
   (0.028)~~ (0.215)\\

   \\
  0.047 +  0.681 Y_{t-1} , \text{ for } t > 313, \\
  (0.013)~~ (0.230)
\end{array}
\right.
\]
where in parentheses are the standard errors of the estimators obtained from the robust sandwich matrix $\widehat H^{-1}_j \widehat \Sigma^{-1}_j \widehat H^{-1}_j$ computed on each regime $j$, where  $\widehat \Sigma_j$ is given by (\ref{eq_Sigma_j}) and
$\widehat H_j =  \frac{1}{n}\sum_{t=1}^{n} \frac{\partial^2 \widehat \ell_{t,j}(\widehat{\theta}_n(\widehat T_j))}{ \partial \theta \partial \theta' } $.
These parameters estimation display a distortion in term of standard errors ; it can be explained by the fact that the true distribution of the observations (which is binary),
is quite different from the Poisson quasi-likelihood used.

\subsection{Number of trades in the stock of Technofirst}

\begin{figure}[h!]
\begin{center}
\includegraphics[height=6cm, width=11.5cm]{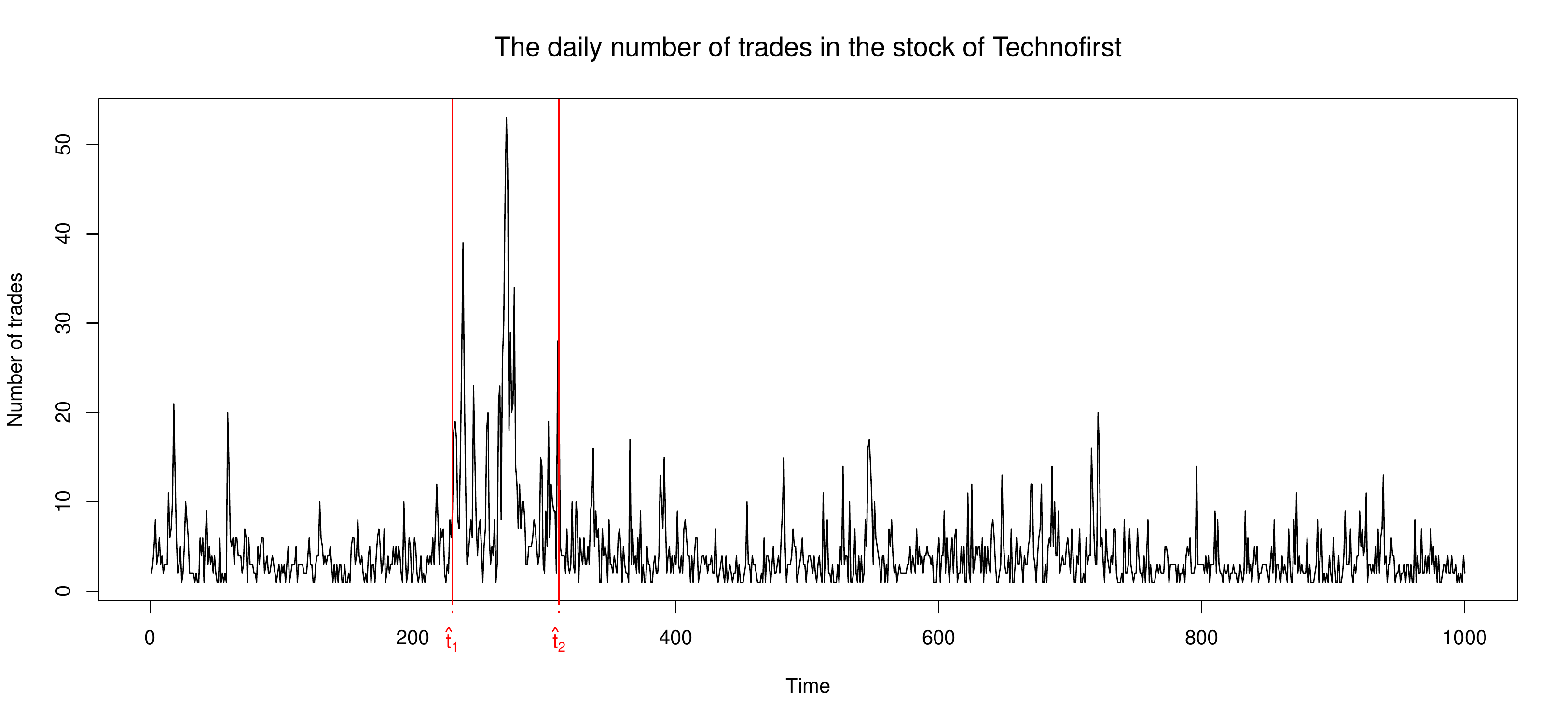}
\end{center}
\vspace{-1cm}
\caption{ The daily number of trades in the stock of the Technofirst with the estimated locations of the breakpoints $\widehat{t}_1$ and $\widehat{t}_2$. }
\label{Result_Technofirst_2}
\end{figure}

 Secondly, we apply our change-point detection procedure to a financial time series data. We consider the daily number of trades in the stock of Technofirst listed in the
 NYSE Euronext group.
 It is a series of $1000$ observations  from \emph{04 January 2010} to \emph{20 April 2016} (see Figure \ref{Result_Technofirst_2}).
 The data are available online at the website "https://www.euronext.com/en/products/equities/FR0011651819-ALXP".
These data have been analyzed by Ahmad and Francq \cite{Francq2016} with the PQMLE, and have concluded that the INGARCH(1,3) is more appropriate.
We carry out an INGARCH(1,1) with the possibility of change in the observations.

\medskip

\begin{figure}[h!]
\begin{center}
\includegraphics[height=5.5cm, width=15.cm]{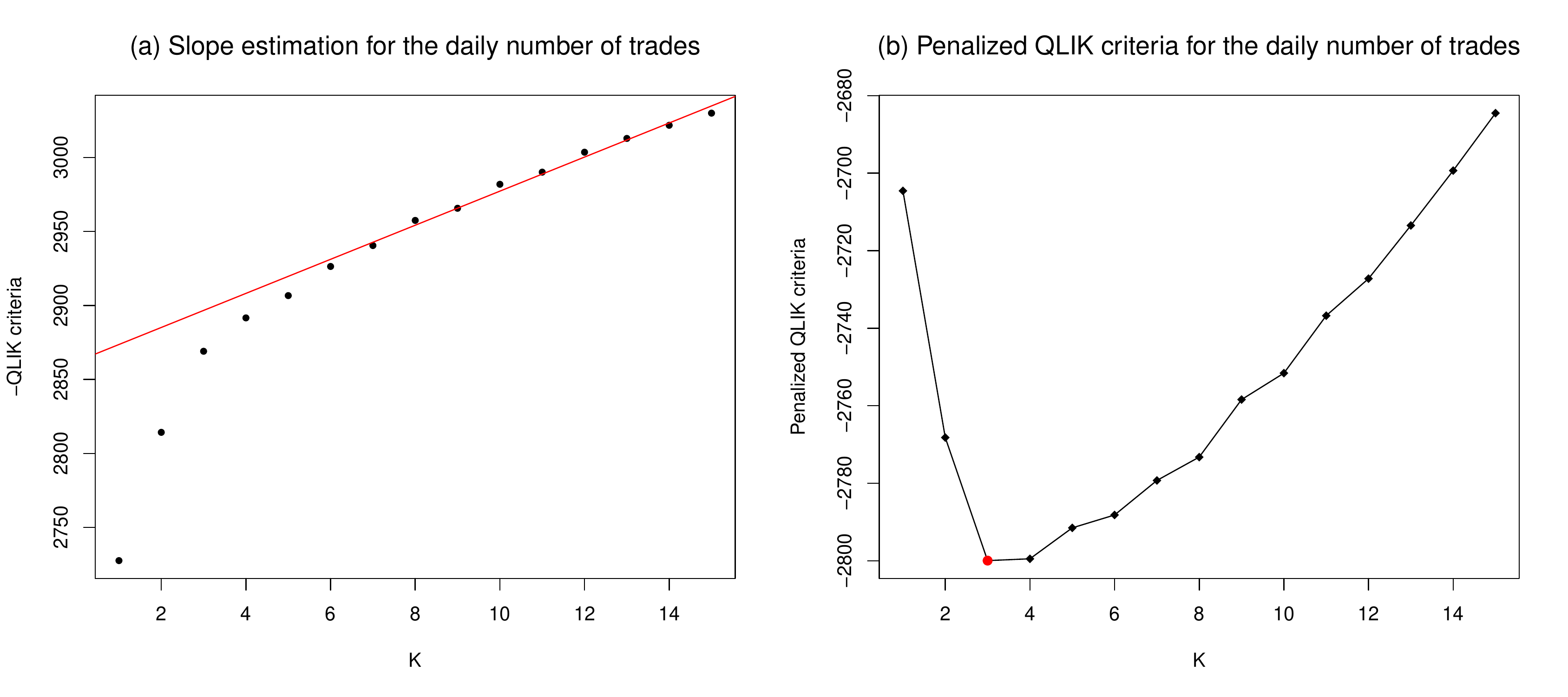}
\end{center}
\vspace{-1cm}
\caption{ The curve of $-\min_{\underline{t},\underline{\theta}} QLIK(K)$ and the graph $(K,\min_{\underline{t},\underline{\theta}} penQLIK(K))$ for the daily number of trades in the stock of Technofirst.
The solid line represents the linear part of this curve with slope $\widehat{\kappa}_n/2 = 11.52$.}
\label{Result_Technofirst_1}
\end{figure}

\noindent 
The slope estimation procedure obtained with $u_n=\left[(\log (n))^{2} \right]$ returns  $\widehat{\kappa}_n \approx 23.04$ and the estimation of the number of segments is
$\widehat{K}_n=3$, {\it i.e.} two changes are detected (see Figure \ref{Result_Technofirst_1}).
The locations of the breakpoints estimated are $\widehat{t}_1=230$ (06 April 2011) and $\widehat{t}_2=311$ (06 September 2011),
see also Figure \ref{Result_Technofirst_2}.
%
%

\medskip

\noindent The estimated model with  change-points is
\begin{equation*}
\E(Y_t | \mathcal{F}_{t-1})=\lambda_{t}=\left\{
\begin{array}{l}
  2.436 + 0.368 Y_{t-1} , \text{ for } t \leq 230 \\
   (0.126)~~ (0.032) \\
   \\
   4.643 + 0.607 Y_{t-1} +0.032\lambda_{t-1}, \text{ for } 230 < t \leq  311\\
  (0.649)~~ (0.022)~~~~~ (0.033)\\
   \\
  1.113 + 0.166 Y_{t-1} +0.531 \lambda_{t-1}, \text{ for } t > 311, \\
  (0.226)~~ (0.016)~~~~~(0.071)
\end{array}
\right.
\end{equation*}
where in parentheses are the robust standard errors of the estimators.
Let us note that, we have applied the TNOC,  which fund that the INARCH(1) representation is the most appropriate for first regime.

\section{Proofs of the main results}


 In the sequel, $C$ denotes a positive constant whom value may differ from an inequality to another
 and we set $v_n=n/\kappa_n$ for all $n\geq 1$.

\subsection{Proof of Proposition \ref{prop1}}
Remark that since $\{Y_{t,j},t\in \Z\}$ is stationary and ergodic, the process $\{\ell_{t,j}(\theta),t\in \Z\}$ is also a stationary and ergodic sequence,
for any $\theta \in \Theta$. Then, the proof can be  divided into two parts. For $j=1,\cdots,K^*$, we will first {\bf (1.)} show that
$\big\| \frac{1}{n^*_j}\widehat{L}_{n}(T^*_j, \theta)-\E (\ell_{1,j}(\theta)) \big\|_\Theta \limitepsn 0 $ ;
secondly {\bf (2.)}, we will show that the function $\theta \mapsto \E(\ell_{1,j}(\theta))$ has a
unique maximum in $\theta^*_j$.
\begin{enumerate}
    \item [\bf (1.)] Recall that  $\ell_{t,j}(\theta) = Y_{t,j}\log f^\theta_{t,j}- f^\theta_{t,j}$, for any $\theta \in \Theta$.  We have
\begin{align*}
|\ell_{t,j} | &\leq Y_{t,j} |\log f^\theta_{t,j} | +  |f^\theta_{t,j} |\\
&\leq
 Y_{t,j} \Big|\log \Big( \frac{f^\theta_{t,j}}{\underline{c}} \times \underline{c} \Big) \Big| +  |f^\theta_{t,j} |\\
&\leq
 Y_{t,j} \big(\Big|\frac{f^\theta_{t,j}}{\underline{c}} -1\Big|  + | \log \underline{c} | \big) +  |f^\theta_{t,j} |~~(\textrm{because~for}~x>1,~ |\log x | \leq  | x-1 |)\\
 &\leq
 Y_{t,j} \big(\big|\frac{f^\theta_{t,j}}{\underline{c}} \big| + 1  + | \log \underline{c} | \big) +  |f^\theta_{t,j} |.
\end{align*}
Hence,
\begin{equation}\label{Maojor_N_ell_tj}
 \| \ell_{t,j} \|_\Theta \leq
 Y_{t,j} \big(\frac{1}{\underline{c}} \|f^\theta_{t,j} \|_\Theta +1  + | \log \underline{c} | \big) +  \|f^\theta_{t,j}  \|_\Theta .
\end{equation}
We will show that $\E [ \| \ell_{t,j} \|_\Theta ] < \infty$. According to (\ref{Maojor_N_ell_tj}), we have
\begin{align*}
\E [ \| \ell_{t,j}  \|_\Theta ] & \leq
\E \big[Y_{t,j} \big(\frac{1}{\underline{c}}  \|f^\theta_{t,j} \|_\Theta +1  + | \log \underline{c} | \big) +  \|f^\theta_{t,j} \|_\Theta  \big]
\leq
C \E \big[ \big( \frac{Y_{t,j}}{\underline{c}} +1 \big) \|f^\theta_{t,j} \|_\Theta  \big]\\
&\leq
  C\big(\E \big[ \big(\frac{Y_{t,j}}{\underline{c}}  +1 \big)^2 \big] \big)^{1/2} (\E \| f^\theta_{t,j} \|^2_\Theta )^{1/2} \leq
  C (\E  \| f^\theta_{t,j} \|^2_\Theta )^{1/2}.
\end{align*}
Since (\textbf{A$_0(\Theta)$}) holds, we have
\begin{align*}
 \|f^\theta_{t,j} \|_\Theta & \leq  \|f^\theta_{t,j} - f^\theta(0,\cdots,0) |_\Theta + \|f^\theta(0,\cdots,0) \|_\Theta
 \leq  \sum\limits_{k \geq 1} \alpha^{(0)}_k |Y_{t-k,j}|  + \|f^\theta(0,\cdots,0) \|_\Theta.
\end{align*}
Therefore,
\begin{align*}
\E [ \| \ell_{t,j} \|_\Theta ] & \leq
C \big[\E \big( \sum\limits_{k \geq 1} \alpha^{(0)}_k |Y_{t-k,j}|  + \|f^\theta(0,\cdots,0) \|_\Theta \big)^2 \big]^{1/2} \\
&\leq
C   \sum\limits_{k \geq 1} \alpha^{(0)}_k (\E |Y_{t-k,j}|^2)^{1/2}  +(\E \|f^\theta(0,\cdots,0) \|^2_\Theta )^{1/2}\\
&\leq
C   \sum\limits_{k \geq 1} \alpha^{(0)}_k +(\E \|f^\theta(0,\cdots,0) \|^2_\Theta )^{1/2} <\infty.
\end{align*}
By the uniform strong law of large number applied on the process  $\{\ell_{t,j}(\theta),t\in \Z\}$, it holds that
\begin{equation} \label{conv2}
\big\| \frac{1}{n^*_j} L_{n,j}(T^*_j, \theta)-\E (\ell_{1,j}(\theta)) \big\|_\Theta= \big\| \frac{1}{n^*_j} \sum_{t\in T^*_j} \ell_{t,j}(\theta)-\E (\ell_{1,j}(\theta)) \big\|_\Theta  \limitepsn 0 .
 \end{equation}
 According to Ahmad and Francq \cite{Francq2016}, we have
\begin{equation}\label{conv1_Francq}
 \frac{1}{n^*_j} \big\| \widehat{L}_{n,j}(T^*_j, \theta)-L_{n,j}(T^*_j, \theta) \big\|_\Theta \limitepsn 0,
 \text{ for any } j = 1,\cdots,K^*.
 \end{equation}
From (\ref{conv2}) and (\ref{conv1_Francq}), we deduce that
\begin{equation} \label{conv3}
\big\| \frac{1}{n^*_j} \widehat{L}_{n,j}(T^*_j, \theta)-\E (\ell_{1,j}(\theta)) \big\|_\Theta \limitepsn 0 .
 \end{equation}
The following lemma is needed to complete the proof of {\bf (1.)}.
\begin{lem} \label{lem0}
 Under the assumptions of Theorem \ref{th1}, for any $j=1,\cdots,K^*$,
\[
 \frac{ v_{n^*_j} }{n^*_j} \big\| \widehat{L}_n(T^*_j, \theta)-\widehat{L}_{n,j}(T^*_j, \theta) \big\|_\Theta \limitepsn 0 .
 \]
\end{lem}
 From  (\ref{conv3}) and Lemma \ref{lem0}, we have

 \begin{equation} \label{conv4}
\big \| \frac{1}{n^*_j} \widehat{L}_{n}(T^*_j, \theta)-\E (\ell_{1,j}(\theta)) \big\|_\Theta \limitepsn 0 .
 \end{equation}
 \item [\bf (2.)] Now, we show that the  function $\theta \mapsto \E(\ell_{1,j}(\theta))$ has a unique maximum in $\theta^*_j$.
 We will proceed as in Doukhan and Kengne (2015).
 For any $\theta \in \Theta$, define  $L^{(j)}(\theta) :=\E[\ell_{1,j}(\theta)]$. Let $\theta \in \Theta$, with $\theta \neq \theta^*_j$. We have
 \begin{align*}
L^{(j)}(\theta^*_j) - L^{(j)}(\theta) &= \E[\ell_{1,j}(\theta^*_j)] -\E[\ell_{1,j}(\theta)] \\
&= \E[Y_{1,j}\log f^{\theta ^*_j}_{1,j}- f^{\theta ^*_j}_{1,j}] -\E[Y_{1,j}\log f^\theta_{1,j}- f^\theta_{1,j}]\\
&= \E [ f^{\theta ^*_j}_{1,j} (\log f^{\theta ^*_j}_{1,j} -\log f^\theta_{1,j} ) ] - \E [f^{\theta ^*_j}_{1,j} - f^\theta_{1,j}].
\end{align*}
We apply the mean value theorem at the function $x \mapsto \log x$ defined in $[\underline{c},+\infty[$.
There exists $\xi$ between $f^{\theta ^*_j}_{1,j}$ and $f^{\theta}_{1,j}$ such that
 \[
 \log f^{\theta ^*_j}_{1,j} -\log f^\theta_{1,j}= \frac{1}{\xi}(f^{\theta ^*_j}_{1,j} - f^\theta_{1,j}).
 \]
Hence,
 \begin{align*}
L^{(j)}(\theta^*_j) - L^{(j)}(\theta)
&= \E \Big[ \frac{f^{\theta ^*_j}_{1,j}}{\xi} (f^{\theta ^*_j}_{1,j} - f^\theta_{1,j} ) \Big] - \E \Big[f^{\theta ^*_j}_{1,j} - f^\theta_{1,j} \Big]\\
&= \E \Big[ \Big( \frac{f^{\theta ^*_j}_{1,j}}{\xi} -1 \Big)\left(f^{\theta ^*_j}_{1,j} - f^\theta_{1,j}\right) \Big]
= \E \big[ \frac{1}{\xi} ( f^{\theta ^*_j}_{1,j}-\xi )(f^{\theta ^*_j}_{1,j} - f^\theta_{1,j} ) \big].
\end{align*}
From assumption (\textbf{A0}), it follows that  $\frac{1}{\xi}( f^{\theta ^*_j}_{1,j}-\xi  )(f^{\theta ^*_j}_{1,j} - f^\theta_{1,j}) \neq 0 $ a.s,  since $\theta \neq \theta^*_j$.\\
Moreover,
\begin{itemize}
    \item if $f^{\theta ^*_j}_{1,j} <f^{\theta}_{1,j}$, then $f^{\theta ^*_j}_{1,j}< \xi <f^{\theta}_{1,j}$ and hence $\frac{1}{\xi}( f^{\theta ^*_j}_{1,j}-\xi  )(f^{\theta ^*_j}_{1,j} - f^\theta_{1,j}) > 0 $ ;

    \item if $f^{\theta ^*_j}_{1,j} >f^{\theta}_{1,j}$, then $f^{\theta}_{1,j}< \xi <f^{\theta^*_j}_{1,j}$ and hence $\frac{1}{\xi}( f^{\theta ^*_j}_{1,j}-\xi  )(f^{\theta ^*_j}_{1,j} - f^\theta_{1,j}) > 0 $.
\end{itemize}
We deduce that  $\frac{1}{\xi}( f^{\theta ^*_j}_{1,j}-\xi  )(f^{\theta ^*_j}_{1,j} - f^\theta_{1,j}) > 0 $ a.s. Hence, $L^{(j)}(\theta^*_j) - L^{(j)}(\theta)>0$ and the  function $\theta \mapsto \E(\ell_{1,j}(\theta))$ has a unique maximum in $\theta^*_j$.\\
~\\
{\bf (1.)}, {\bf (2.)} and the standard arguments lead to the conclusion. $~~~~~~~~~~~~~~~~~~~~~~~~~~~~~~~~~~~~~~~~~~~~~~~~~~~~~~~~~~~~~~~~ \blacksquare$
\end{enumerate}
~~\\
\emph{\bf Proof of Lemma \ref{lem0}}\\
For any $j=1,\cdots,K^*$, remark that
\begin{align*}
  \| \widehat{L}_n(T^*_j, \theta)-\widehat{L}_{n,j}(T^*_j, \theta)\|_\Theta
 &\leq
 \sum_{t \in T^*_j}\|\widehat{\ell}_{t}(\theta)-\widehat{\ell}_{t,j}(\theta) \|_\Theta\\
 &\leq
  \sum_{t \in T^*_j} \|Y_{t,j}\log \widehat{\lambda}_{t}(\theta)- \widehat{\lambda}_{t}(\theta) -Y_{t,j}\log \widehat{\lambda}_{t,j}(\theta)+ \widehat{\lambda}_{t,j}(\theta)\|_\Theta\\
  &\leq
  \sum_{t \in T^*_j}  (Y_{t,j} \|\log \widehat{f}^\theta_t-\log \widehat{f}^\theta_{t,j}\|_\Theta +\| \widehat{f}^\theta_t- \widehat{f}^\theta_{t,j} \|_\Theta ).
 \end{align*}
According to the proprieties of the function $x \mapsto \log x$, we can show that
 $ \| \log \widehat{f}^\theta_t-\log \widehat{f}^\theta_{t,j}\|_\Theta \leq \frac{1}{\underline{c}}\| \widehat{f}^\theta_t- \widehat{f}^\theta_{t,j} \|_\Theta$.
Moreover, for $ t \in T^*_j$, we have
\begin{align}
  \| \widehat{f}^\theta_t- \widehat{f}^\theta_{t,j} \|_\Theta
   \nonumber   &=  \| f ( Y_{t-1},\cdots,Y_{1},0,\cdots;\theta) - f ( Y_{t-1,j},\cdots,Y_{t^*_{j-1}+1,j},0,\cdots ;\theta) \|_\Theta \\
     \label{proof_lem0_eq1} & \leq \sum\limits_{k=t-t^*_{j-1}}^{t-1} \alpha^{(0)}_{k} |Y_{t-k}| \leq  \sum\limits_{k \geq t-t^*_{j-1} } \alpha^{(0)}_{k} Y_{t-k} .
\end{align}
 Hence
 \begin{align*}
 \frac{v_{n^*_j}}{n^*_j} \| \widehat{L}_n(T^*_j, \theta)-\widehat{L}_{n,j}(T^*_j, \theta) \|_\Theta
 & \leq   \frac{v_{n^*_j}}{n^*_{j}} \sum_{t \in T^*_j}   \Big[ \Big( \frac{Y_{t,j}}{\underline{c}}+1 \Big) \| \widehat{f}^\theta_t- \widehat{f}^\theta_{t,j} \|_\Theta \Big] \\
 &\leq   \frac{v_{n^*_j}}{n^*_{j}}  \sum_{t \in T^*_j}  \sum\limits_{k \geq t-t^*_{j-1} } \alpha^{(0)}_{k}  \Big[ \Big( \frac{Y_{t,j}}{\underline{c}}+1 \Big) Y_{t-k} \Big] \\
 &\leq   \frac{v_{n^*_j}}{n^*_{j}} \sum\limits_{ \ell=1 }^{n^*_j}  \sum\limits_{k \geq \ell } \alpha^{(0)}_{k}  \Big[ \Big( \frac{ Y_{\ell + t^*_{j-1},j} }{\underline{c}}+1 \Big) Y_{\ell + t^*_{j-1} -k} \Big].
 \end{align*}
By using  Kounias and Weng (1969), it suffices to show that
\[ \sum_{\ell \geq 1}\frac{v_{\ell} }{ \ell } \E \Big[ \sum\limits_{k \geq \ell } \alpha^{(0)}_{k} \Big( \frac{ Y_{\ell + t^*_{j-1},j} }{\underline{c}}+1\Big) Y_{\ell + t^*_{j-1} -k}   \Big]  < \infty.  \]

By using Hölder's inequality, for any $\ell \geq 1$, $k \geq \ell$, it holds that (see (\ref{moment}))
 \[ \E  \Big[ \Big( \frac{ Y_{\ell + t^*_{j-1},j} }{\underline{c}}+1 \Big) Y_{\ell + t^*_{j-1} -k} \Big]  \leq
  \Big( \E\Big[ \Big( \frac{ Y_{\ell + t^*_{j-1},j} }{\underline{c}}+1 \Big)^2 \Big]  \Big)^{1/2}   \times  (\E Y^2_{\ell + t^*_{j-1} -k} )^{1/2}
  = C < \infty .\]
Hence,
\[
 \sum_{\ell \geq 1}\frac{v_{\ell} }{ \ell } \E \Big[ \sum\limits_{k \geq \ell } \alpha^{(0)}_{k} \Big( \frac{ Y_{\ell + t^*_{j-1},j} }{\underline{c}}+1\Big) Y_{\ell + t^*_{j-1} -k}   \Big]
 \leq C \sum_{\ell \geq 1}  \frac{1 }{ \kappa_{\ell }}   \sum\limits_{k \geq \ell }   \alpha^{(0)}_{k} < \infty  ,
  \]
  where the last equation follows from  assumption (\ref{eq_th1}) on the regularization parameter.
  Thus,
  \[ \frac{v_{n^*_j}}{n^*_j}\big\| \widehat{L}_n(T^*_j, \theta)-\widehat{L}_{n,j}(T^*_j, \theta)\big\|_\Theta \limitepsn 0 . \]
  $~~~~~~~~~~~~~~~~~~~~~~~~~~~~~~~~~~~~~~~~~~~~~~~~~~~~~~~~~~~~~~~~~~~~~~~~~~~~~~~~~~~~~~~~~~~~~~~~~~~~~~~~~~~~~~~~~~~~~~~~~~~~~~~~~~~~~~~~~~~~~~~~ \blacksquare$\\


\subsection{Proof of Theorem \ref{th1}}
We will proceed as in Bardet \emph{et al.} (2012). Firstly, we assume that $K^*$ is known
and we show $(\widehat{\underline{\tau}}_n,\widehat{\underline{\theta}}_n)  \limiteproban  \left(\underline{\tau}^*,\underline{\theta} ^*\right)$.
Secondly, $K^*$ is assumed to be unknown and we show $\widehat{K}_n  \limiteproban  K^*$; which completes the proof of the theorem.
\begin{enumerate}
\item [\bf (1.)] Assume that $K^*$ is known and denote for any $\underline{t} \in \mathcal{M}_n(K^*)$:
    \[\widehat{I}_n(\underline{t}):=\widehat{J}_n(K^*,\underline{t},\widehat{\underline{\theta}}_n(\underline{t})) =-2 \sum_{k=1}^{K^*} \sum_{j=1}^{K^*}  \widehat L_n ( T_k \cap T^*_j,\widehat \theta_n(T_k) ). \]
Hence,
  $\widehat{\underline{t}}_n =  \underset{ \underline{t} \in \mathcal{M}_n(K^*)}{\text{argmin}} ( \widehat{I}_n(\underline{t})) $.
 Let us show that $\widehat{\underline{\tau}}_n \limiteproban \underline{\tau}^*$; which will implies that
$\widehat{\theta}_n (\widehat T_{n,j}) \limiteproban \widehat{\theta}_n (T^*_{j})$ and from Proposition \ref{prop1}, $\widehat{\theta}_n (\widehat T_{n,j}) \limiteproban \theta^*_{j}$ for all
$j = 1, \cdots ,K^*$.
 Without loss of generality, we assume that $K^* = 2$.  Let $t^*$ be the change-point location and $(u_n)_{n \geq 1}$ be a sequence of positive integers satisfying
 $u_n \rightarrow \infty ,~u_n/n \rightarrow 0$. For some $0<\eta<1$, define
\begin{align*}
V_{\eta,u_n}&= \{  t \in \Z ~/ ~| t-t^*|> \eta n~;~u_n \leq t \leq n-u_n\},\\
W_{\eta,u_n}&= \{  t \in \Z ~/ ~| t-t^*|> \eta n~;~0 < t <u_n~\textrm{or}~ n-u_n < t \leq n\}.
\end{align*}
Remark that we have asymptotically $\prob ( \| \widehat{\underline{\tau}}_n - \underline{\tau}^* \|_m    > \eta) \simeq \prob (\|\widehat{t}_n - t^* \|_m    > \eta n ).  $
 But
 \begin{align*}
 \prob (\|\widehat{t}_n - t^* \|_m    > \eta n )
 & \leq \prob (\widehat{t}_n \in V_{\eta,u_n} ) + \prob (\widehat{t}_n \in W_{\eta,u_n} )\\
 & \leq \prob ( \underset{ t \in V_{\eta,u_n} }{\min} ( \widehat{I}_n(t) - \widehat{I}_n(t^*)) \leq 0 ) +
 \prob ( \underset{ t \in W_{\eta,u_n} }{\min} ( \widehat{I}_n(t) - \widehat{I}_n(t^*)) \leq 0 ).
 \end{align*}
We will show that these two probabilities tend to 0.
Let us show that $\prob ( \underset{ t \in V_{\eta,u_n} }{\min} ( \widehat{I}_n(t) - \widehat{I}_n(t^*)) \leq 0 ) \rightarrow 0$.
Let $t \in V_{\eta,u_n}$ satisfying $t \geq t^*$  (without loss of generality). Then, we have $T_1 \cap T^*_1=T^*_1$, $T_2 \cap T^*_1=\emptyset$ and $T_2 \cap T^*_2=T_2$.
 We have the decomposition
 \begin{align}\label{Dif_I}
 \widehat{I}_n(t) - \widehat{I}_n(t^*)
  & =
 2 \big[
 \widehat L_n (T^*_1 ,\widehat \theta_n(T^*_1) ) -\widehat L_n (T^*_1 ,\widehat \theta_n(T_1) )
  -\widehat{L}_n ( T_1 \cap T^*_2,\widehat \theta_n(T_1) ) \nonumber \\
  &~~~~~~~~+
   \widehat L_n (T_2,\widehat \theta_n(T^*_2)) -
   \widehat L_n (T_2,\widehat \theta_n(T_2)) + \widehat{L}_n ( T_1 \cap T^*_2,\widehat \theta_n(T^*_2) )
   \big].
 \end{align}
Since $| T^*_1 |=t^*$, $ | T_1 \cap T^*_2 | =t-t^*$ and  $ |  T_2 | = n-t \geq u_n$,  each term tends to $\infty$ with $n$.\\
Recall that for $j=1,2$,  $L^{(j)}(\theta)=\E(\ell_{1,j}(\theta))$.
According to Proposition \ref{prop1} and relation (\ref{conv4}), we get the following convergences, uniformly on $V_{\eta,u_n}$,
 \begin{align*}
 ~~& \widehat{\theta}_n(T^*_1) \limitepsn \theta^*_1,~~\widehat{\theta}_n(T^*_2) \limitepsn \theta^*_2,
 ~~ \widehat{\theta}_n(T_2) \limitepsn \theta^*_2,~~~
 \big\| \frac{1}{n} \widehat{L}_{n}(T^*_1, \theta)-\tau^*_1 L^{(1)}(\theta) \big\|_\Theta \limitepsn 0, ~~\\
~&
 \big\| \frac{1}{t-t^*} \widehat{L}_{n}(T_1 \cap T^*_2, \theta)- L^{(2)}(\theta) \big \|_\Theta \limitepsn 0 \text{ and }
 ~~
\big\| \frac{1}{n-t} \widehat{L}_{n}(T_2, \theta)- L^{(2)}(\theta)  \big \|_\Theta \limitepsn 0 .
 \end{align*}
 Set $\eta_n = (t-t^*)/n $ ; clearly, $\varepsilon_n \in (\eta, 1)$. From (\ref{Dif_I}), we get
 \begin{align}\label{Dif_I_conv}
 \nonumber \frac{1}{2n}(\widehat{I}_n(t) - \widehat{I}_n(t^*) ) &=
 \frac{1}{n}\big( \widehat L_n (T^*_1 ,\widehat \theta_n(T^*_1) ) -\widehat L_n (T^*_1 ,\widehat \theta_n(T_1) )  \big)
  + \frac{1}{n}\big( \widehat L_n (T_2,\widehat \theta_n(T^*_2)) - \widehat L_n (T_2,\widehat \theta_n(T_2)) \big) \\
\nonumber   & \hspace{4.5cm} + \frac{1}{n}\big( \widehat{L}_n ( T_1 \cap T^*_2,\widehat \theta_n(T^*_2) ) - \widehat{L}_n ( T_1 \cap T^*_2,\widehat \theta_n(T_1) ) \big) \\
\nonumber & = \tau_1^*  \big( L^{(1)}(\theta_1^*) - L^{(1)}(\widehat \theta_n(T_1) )  \big) + o(1) + o(1)
   + \eta_n \big( L^{(2)}(\theta_2^*)  - L^{(2)} (\widehat \theta_n(T_1) ) \big) + o(1) \\
 &= \tau_1^*  \big( L^{(1)}(\theta_1^*) - L^{(1)}(\widehat \theta_n(T_1) )  \big)
   + \eta_n \big( L^{(2)}(\theta_2^*)  - L^{(2)} (\widehat \theta_n(T_1) ) \big) + o(1)  .
 \end{align}
 Let $\mathcal{V}_1$ and $\mathcal{V}_2$ be two disjoint open neighborhoods of $\theta^*_1$ and $\theta^*_2$ respectively.  For $j=1,2$, define
 \[ \delta_j= \underset{\theta \in  \mathcal{V}^c_j}{\inf} \big( L^{(j)}(\theta^*_j) - L^{(j)}(\theta) \big). \]
 Remark that $\delta_j >0$ since the function $\theta \mapsto  L^{(j)}(\theta) $  has a strict maximum in $\theta^*_j$ (see proof of Proposition \ref{prop1}).
  With $\varepsilon_n=\min(\tau^*_1\delta_1,\eta_n \delta_2)$ and $\varepsilon=\min(\tau^*_1\delta_1,\eta \delta_2)$, we have
\begin{itemize}
    \item if $\widehat{\theta}_n(T_1) \in \mathcal{V}_1$, that is $\widehat{\theta}_n(T_1) \in \mathcal{V}^c_2$,
    then $  \eta_n \big( L^{(2)}(\theta_2^*)  - L^{(2)} (\widehat \theta_n(T_1) ) \big) >\eta_n \delta_2 $;
    \item if $\widehat{\theta}_n(T_1) \notin \mathcal{V}_1$, that is $\widehat{\theta}_n(T_1) \in \mathcal{V}^c_1$,
    then $\tau_1^*  \big( L^{(1)}(\theta_1^*) - L^{(1)}(\widehat \theta_n(T_1) )  \big) >\tau^*_1\delta_1$.
\end{itemize}
 In both cases,  $ \frac{1}{2n} (\widehat{I}_n(t) - \widehat{I}_n(t^*) ) \geq \varepsilon_n + o(1) \geq \varepsilon + o(1)$, for any $t \in V_{\eta,u_n}$.
 This implies $\prob \Big( \underset{ t \in V_{\eta,u_n} }{\min} \big( \widehat{I}_n(t) - \widehat{I}_n(t^*)\big) \leq 0 \Big) \rightarrow 0$.
    By going along similar lines, one can prove that $\prob \Big( \underset{ t \in W_{\eta,u_n} }{\min} \big( \widehat{I}_n(t) - \widehat{I}_n(t^*)\big) \leq 0 \Big) \rightarrow 0$.
    Then, it follows that $\eta >0$, $\prob \big(  \| \widehat{\underline{\tau}}_n - \underline{\tau}^* \|_m    > \eta \big) \rightarrow 0  ~~ \textrm{as}~~n\rightarrow \infty$.\\
 \item [\bf (2.)] Now, assume that $K^*$ is  unknown. For $K \geq 2$, $x = (x_1, \cdots , x_{K-1}) \in \R^{K-1}$,
 $y = (y_1, \cdots , y_{K^*-1}) \in \R^{K^*-1}$, denote
 \[ \left\| x-y\right\|_\infty =\underset{1 \leq j \leq K^*-1}{\max} ~~\underset{1 \leq k \leq K-1}{\min} |x_k-y_j|. \]
The following lemma will be useful in the sequel. It follows from {\bf (1.)} and the definition of $\left\| \cdot \right\|_\infty$.
\begin{lem}\label{lem1}
 Let $K \geq 1$, $( \widehat{\underline{t}}_n,  \widehat{\underline{\theta}}_n)$ the estimator obtained by minimizing
 $\widehat{J}_n(K,\underline{t},\underline{\theta})$ on $ \mathcal{M}_n(K) \times \Theta^K$ and $\underline{\widehat{\tau}}_n=\underline{\widehat{t}}_n/n$.
 Under the assumptions of Theorem \ref{th1}, if $K \geq K^*$, then
$\left\| \underline{\widehat{\tau}}_n -\underline{\tau}^* \right\|_\infty \limiteproban  0$.
\end{lem}
We will also use the following lemma, which the proof follows from the Lemma 3.3 of Lavielle and Ludena (2000)
and the argument given in the proof of Lemma 6.5 in \cite{Bardet2012}.
\begin{lem}\label{lem2}
Under the assumptions of Theorem \ref{th1}, for any $K \geq 2$, there
exists $C_K > 0$ such that:
\[
\forall (\underline{t},\underline{\theta}) \in \mathcal{M}_n(K) \times \Theta^K,~~e_n(\underline{t},\underline{\theta}) =
2 \sum_{j=1}^{K^*} \sum_{k=1}^{K} \frac{n_{k,j}}{n} \big(L^{(j)}(\theta^*_j) - L^{(j)}(\theta_k) \big)
\geq \frac{C_K}{n} \| \underline{t} -\underline{t}^* \|_\infty
\]
where $L^{(j)}(\theta) =\E(\ell_{1,j}(\theta))$ for all $\theta \in \Theta$ and $j = 1,\cdots, K^*$.
\end{lem}

 \medskip

Now, let us use the Lemma \ref{lem1}  and \ref{lem2}  to show that $\widehat{K}_n  \limiteproban  K^*$.
To this end, we will show that $\prob (\widehat{K}_n  = K ) \underset{n \rightarrow \infty}{\rightarrow} 0 $, for $K <K^*$ and $ K^* < K \leq K_{\max}$ separately.
In any case, we have
\begin{align}\label{Major1_prob_Kn}
\prob (\widehat{K}_n  = K )
& \leq \prob \Big( \underset{(\underline{t},\underline{\theta}) \in \mathcal{M}_n(K) \times \Theta^K}{\inf} \big( \widetilde{J}_n(K,\underline{t},\underline{\theta})\big) \leq \widetilde{J}_n(K^*,\underline{t}^*,\underline{\theta}^*) \Big)
\nonumber\\
 & \leq
 \prob \Big( \underset{(\underline{t},\underline{\theta}) \in \mathcal{M}_n(K) \times \Theta^K}{\inf} \big( \widehat{J}_n(K,\underline{t},\underline{\theta}) -\widehat{J}_n(K^*,\underline{t}^*,\underline{\theta}^*)\big) \leq \frac{n}{v_n} ( K^*-K) \Big).
\end{align}
\begin{enumerate}
    \item [\bf i-)]  For $K <K^*$, decompose
 $\widehat{J}_n(K,\underline{t},\underline{\theta}) -\widehat{J}_n(K^*,\underline{t}^*,\underline{\theta}^*)=n \left(d_n (\underline{t},\underline{\theta})+ e_n(\underline{t},\underline{\theta})\right)$,
 where $e_n(\underline{t},\underline{\theta})$ is defined in Lemma \ref{lem2}  and
\[
d_n (\underline{t},\underline{\theta})
= 2 \Big[
\sum_{j=1}^{K^*} \frac{n^*_{j}}{n} \big(\frac{\widehat L_n \left( T^*_j, \theta^*_j\right)}{n^*_{j}}- L^{(j)}(\theta^*_j) \big)
+ \sum_{k=1}^{K}\sum_{j=1}^{K^*} \frac{n_{k,j}}{n} \big(L^{(j)}(\theta_k) - \frac{\widehat L_n \left( T^*_j \cap T_k, \theta_k\right)}{n_{k,j}} \big) \Big].
\]
Hence, from (\ref{Major1_prob_Kn}), we have
\begin{equation}\label{Major2_prob_Kn}
\prob (\widehat{K}_n  = K ) \leq
\prob \big( \underset{(\underline{t},\underline{\theta}) \in \mathcal{M}_n(K) \times \Theta^K}{\inf} (d_n (\underline{t},\underline{\theta})+ e_n(\underline{t},\underline{\theta}) ) \leq \frac{1}{v_n}  ( K^*-K ) \big).
\end{equation}
The equation (\ref{conv3}) ensures that $d_n (\underline{t},\underline{\theta}) \rightarrow 0 $ a.s. and uniformly on $\mathcal{M}_n(K) \times \Theta^K$.
According to Lemma \ref{lem2}, there exists $C_K>0$ such that $e_n(\underline{t},\underline{\theta}) \geq \frac{C_K}{n} \| \underline{t} -\underline{t}^* \|_\infty$,
for all $(\underline{t},\underline{\theta}) \in \mathcal{M}_n(K) \times \Theta^K$.
Since $K< K^*$, for any $\underline{t} \in \mathcal{M}_n(K)$, we have
$ \frac{1}{n} \| \underline{t} -\underline{t}^*  \|_\infty = \| \underline{\tau} -\underline{\tau}^*  \|_\infty \geq  \underset{1 \leq j \leq K^*}{\min}  ( \tau^*_j-\tau^*_{j-1} )/2 >0$.
Then $e_n(\underline{t},\underline{\theta})>0$ for $(\underline{t},\underline{\theta}) \in \mathcal{M}_n(K) \times \Theta^K$  and since $\frac{1}{v_n} \underset{n \rightarrow \infty}{\longrightarrow} 0 $, we deduce from (\ref{Major2_prob_Kn}) that
$\prob (\widehat{K}_n  = K ) \underset{n \rightarrow \infty}{\longrightarrow} 0$.
\item [\bf ii-)] Now, assume that $ K^* < K \leq K_{\max}$. Denote $\widehat{\underline{t}}_n =(\widehat{t}_{n,1}, \cdots,\widehat{t}_{n,K^*} )$.  From (\ref{Major1_prob_Kn})
and the Markov's inequality, we have:
\begin{align}\label{Major3_prob_Kn}
\prob (\widehat{K}_n  = K )
 & \leq
 \prob \big( \widehat{J}_n(K,\widehat{\underline{t}}_n, \widehat{\underline{\theta}}_n) -\widehat{J}_n(K^*,\underline{t}^*,\underline{\theta}^*) +\frac{n}{v_n}  ( K-K^* ) \leq 0 \big)  \nonumber\\
 & \leq
 \prob \big( | \widehat{J}_n(K,\widehat{\underline{t}}_n, \widehat{\underline{\theta}}_n) -\widehat{J}_n(K^*,\underline{t}^*,\underline{\theta}^*) | > \frac{n}{v_n}  \big)  \nonumber\\
  & \leq
  \frac{v_n}{n}  \E | \widehat{J}_n(K,\widehat{\underline{t}}_n, \widehat{\underline{\theta}}_n) -\widehat{J}_n(K^*,\underline{t}^*,\underline{\theta}^*) |.
\end{align}
By Lemma \ref{lem1}, there exists some subset $\{k_j,~1 \leq j \leq K^*-1\} \subset \{1,\cdots,K-1\}$ such that
for any $j = 1,\cdots,K^*-1$, $\widehat{t}_{n,k_j}/n \rightarrow \tau^*_j$. Set $k_0=0$ and $k_{K^*}=K$. We have
\begin{align*}
\widehat{J}_n(K,\widehat{\underline{t}}_n, \widehat{\underline{\theta}}_n) -\widehat{J}_n(K^*,\underline{t}^*,\underline{\theta}^*)
& = 2 \big( \sum_{j=1}^{K^*}  \widehat L_n(T^*_j,\theta^*_j) -\sum_{k=1}^{K}  \widehat L_n(\widehat T_{n,k},\widehat \theta_{n,k}) \big)\\
&=
2
\sum_{j=1}^{K^*}  \big[\widehat L_n(T^*_j,\theta^*_j) -\sum_{k=k_{j-1} +1}^{k_{j}}  \widehat L_n(\widehat T_{n,k},\widehat \theta_{n,k}) \big]
\end{align*}
and from (\ref{Major3_prob_Kn}), it follows that
\begin{align*}
\prob (\widehat{K}_n  = K ) & \leq
\frac{2v_n}{n}  \sum_{j=1}^{K^*}  \E \big|\widehat L_n(T^*_j,\theta^*_j) -\sum_{k=k_{j-1} +1}^{k_{j}}  \widehat L_n(\widehat T_{n,k},\widehat \theta_{n,k})\big|\\
& \leq
  C \sum_{j=1}^{K^*} \frac{v_{n^*_j}}{n^*_j}\E \big|\widehat L_n(T^*_j,\theta^*_j) -\sum_{k=k_{j-1} +1}^{k_{j}}  \widehat L_n(\widehat T_{n,k},\widehat \theta_{n,k})\big|.
\end{align*}
 For any $j = 1, \cdots ,K^*$, one can easily get from the proof of Lemma \ref{lem0} that
\[
 \frac{v_{n^*_j}}{n^*_j}\E \big|\widehat L_n(T^*_j,\theta^*_j) -\sum_{k=k_{j-1} +1}^{k_{j}}  \widehat L_n(\widehat T_{n,k},\widehat \theta_{n,k})\big| \underset{ n \rightarrow \infty }{\longrightarrow} 0,
 \]
and thus $\prob ( \widehat K_n=K ) \underset{ n \rightarrow \infty }{\longrightarrow} 0$. $~~~~~~~~~~~~~~~~~~~~~~~~~~~~~~~~~~~~~~~~~~~~~~~~~~~~~~~~~~~~~~~~~~~~~~~~~~~~~~~~ \blacksquare$
\end{enumerate}
\end{enumerate}


\subsection{Proof of Theorem \ref{th2}}
 Also for this proof, We will proceed as in \cite{Bardet2012}.
 Without loss of generality, we can assume that $K^* = 2$. Let $(u_n)_{n \geq 1}$ be a sequence satisfying $u_n \limiten \infty$, $\frac{u_n}{n} \limiten 0$
 and $\prob (|\widehat{\underline{t}}_n -\underline{t}^*| > u_n) \limiten 0$ (for instance $u_n= \sqrt{\max(\E |\widehat{\tau}_n - \tau^*|, n^{-1})}$).
 For any $\delta >0$, since we have
\[
\prob (|\widehat{\underline{t}}_n -\underline{t}^*| >\delta)  \leq  \prob (\delta < |\widehat{\underline{t}}_n -\underline{t}^*| \leq u_n) +  \prob (|\widehat{\underline{t}}_n -\underline{t}^*| > u_n),
\]
it suffices to show that
$ \underset{\delta \rightarrow \infty }{\lim}~ \underset{n \rightarrow \infty }{\lim}  \prob (   \delta <  | \widehat{t}_n -t^*| \leq u_n )=0$.

\medskip

Denote $V_{\delta,u_n}= \big\{ t \in \Z / \, \, \delta <  | \widehat{t}_n -t^*| \leq u_n \big\}$. Then
\[
\prob \big( \delta <  | \widehat{t}_n -t^*| \leq u_n \big) \leq \prob \Big(  \underset{t \in V_{\delta,u_n}}{\max} \big(  \widehat{I}_n(t)- \widehat{I}_n(t^*) \big) \leq 0 \Big).
\]
Let $t \in V_{\delta,u_n}$ (for example $t \geq t^*$). With the notation of the proof of Theorem \ref{th1}, we have
$ \widehat{L}_n(T^*_1, \widehat{\theta}_n(T^*_1)) \geq \widehat{L}_n(T^*_1, \widehat{\theta}_n(T_1) )$ and from (\ref{Dif_I}) we get
\[
 \widehat{I}_n(t) - \widehat{I}_n(t^*)  \geq
    \widehat L_n \big( T_1 \cap T^*_2,\widehat \theta_n(T^*_2) \big)  -  \widehat{L}_n \big( T_1 \cap T^*_2,\widehat \theta_n(T_1) \big)
     +  \widehat L_n \big(T_2,\widehat \theta_n(T^*_2) \big) - \widehat L_n \big(T_2,\widehat \theta_n(T_2) \big) .
 \]
 We consider the following two steps.
\begin{enumerate}
    \item [\bf (1.)] Let us show that $ \widehat{L}_n \big( T_1 \cap T^*_2,\widehat \theta_n(T^*_2) \big) -  \widehat{L}_n \big( T_1 \cap T^*_2,\widehat \theta_n(T_1) \big) >0$,
    for $n$ large enough.

 \medskip
   For any $\theta \in \Theta$, we have $\frac{1}{n} \widehat L_n (T_1,\theta )= \frac{t^*}{n} \frac{\widehat L_n (T^*_1,\theta )}{t^*} +
    \frac{t-t^*}{n} \frac{\widehat L_n (T_1 \cap T^*_2,\theta )}{t-t^*}$ and since $ \frac{t-t^*}{n} \leq  \frac{u_n}{n} \limiten 0 $, it follows that
    \[
     \widehat \theta_n(T_1) = \underset{\theta\in \Theta}{\text{argmax}} \big( \frac{1}{n} \widehat L_n (T_1,\theta ) \big) \overset{a.s}{\underset{n,  \delta \rightarrow \infty}{\longrightarrow}} \theta^*_1.
     \]
 Hence, $\frac{1}{t-t^*} \big(\widehat{L}_n ( T_1 \cap T^*_2,\widehat \theta_n(T^*_2) )  -
   \widehat{L}_n ( T_1 \cap T^*_2,\widehat \theta_n(T_1) ) \big)$ converges a.s. and uniformly on $V_{\delta,u_n}$ to $ L^{(2)}(\theta^*_2)- L^{(2)}(\theta^*_1) >0$.

 \medskip
  \item [\bf (2.)] Let us show that  $\frac{1}{t-t^*}  \big(\widehat L_n (T_2,\widehat \theta_n(T^*_2) )  -
  \widehat L_n ( T_2,\widehat \theta_n(T_2) ) \big) \overset{a.s}{\underset{n,  \delta \rightarrow \infty}{\longrightarrow}}   0$.
  For large value of $n$, remark that $\widehat \theta_n(T_2) \in \overset{\circ}{\Theta}$
  so that $\partial  \widehat L_n ( T_2,\widehat \theta_n(T_2) ) / \partial \theta =0$.
   From mean value theorem applied on $\partial  \widehat L_n  / \partial \theta_i$ for any $i=1,\cdots,d$, there exists
   $\widetilde \theta_{n,i} \in [\widehat \theta_n(T_2), \widehat \theta_n(T^*_2)]$ such that
\begin{equation}\label{mvt1}
0=\frac{\partial  \widehat{L}_n ( T_2,\widehat \theta_n(T^*_2) )}{\partial \theta_i} + \frac{\partial^2  \widehat{L}_n ( T_2,\widetilde \theta_{n,i} )}{\partial \theta \partial \theta_i} (\widehat \theta_n(T_2) -\widehat \theta_n(T^*_2) )
\end{equation}
where for $a,b \in \R^d$, $[a,b]=\{ \lambda a + (1-\lambda) b; \, \lambda \in [0,1]\}$.\\
 According the equalities $\widehat{L}_n ( T^*_2, \theta )=\widehat{L}_n ( T_1 \cap T^*_2, \theta ) +
\widehat{L}_n  ( T_2, \theta )$ and $\partial  \widehat{L}_n ( T^*_2,\widehat \theta_n(T^*_2) )/\partial \theta =0$,
it comes from (\ref{mvt1}) that

\[
\frac{\partial  \widehat L_n (T_1 \cap T^*_2,\widehat \theta_n(T^*_2) )}{\partial \theta_i} = \frac{\partial^2  \widehat L_n ( T_2,\widetilde \theta_{n,i} )}{\partial \theta \partial \theta_i} (\widehat \theta_n(T_2) -\widehat \theta_n(T^*_2) ),
 ~~ \text{ for any }  i=1,\cdots,d,
\]
and it follows that
\begin{equation}\label{mvt1.1}
\frac{1}{t-t^*}\frac{\partial \widehat L_n (T_1 \cap T^*_2,\widehat \theta_n(T^*_2) )}{\partial \theta}
= \frac{n-t}{t-t^*} \widetilde{A}_n \cdot (\widehat \theta_n(T_2) -\widehat \theta_n(T^*_2) )
\end{equation}
where $\widetilde{A}_n:=\Big( \frac{1}{n-t} \frac{\partial^2 \widehat L_n ( T_2,\widetilde \theta_{n,i})}{\partial \theta \partial \theta_i} \Big)_{1 \leq i \leq d}$.

\medskip
The following lemma will be useful in the sequel.
\begin{lem}\label{lem3} ~
 %
 %
\begin{itemize}
    \item
     Assume that the conditions of Theorem \ref{th3} hold. Then, for any $j=1,\cdots,K^*$,
%
\begin{equation}\label{lem3_hatL_hatLj}
 \frac{1}{\sqrt{n^*_j}} \Big\| \frac{\partial^i \widehat{L}_n(T^*_j, \theta)}{\partial \theta^i} - \frac{\partial^i \widehat{L}_{n,j}(T^*_j, \theta)}{\partial \theta^i}  \Big\|_\Theta \limitepsn 0 ;
\end{equation}
\item
Assume that the conditions of Theorem \ref{th2} hold. Then, for any $j=1,\cdots,K^*$,

\begin{equation}\label{lem3_hatL_calLj}
\Big\| \frac{1}{n^*_j} \frac{\partial \widehat{L}_n(T^*_j, \theta)}{\partial \theta}-\frac{\partial L^{(j)}(\theta)}{\partial \theta} \Big\|_\Theta \limitepsn 0,
~~ where ~~\frac{\partial L^{(j)}(\theta)}{\partial \theta}=\E\Big(\frac{\partial \ell_{1,j}(\theta)}{\partial \theta} \Big).
\end{equation}
\end{itemize}
\end{lem}

 \medskip

  \medskip

   \medskip

  \medskip

Hence, (\ref{lem3_hatL_calLj}) gives
$
\frac{1}{t-t^*}\frac{\partial  \widehat{L}_n (T_1 \cap T^*_2,\widehat \theta_n(T^*_2) )}{\partial \theta} \overset{a.s}{\underset{n,  \delta \rightarrow \infty}{\longrightarrow}} \frac{\partial L^{(2)}(\theta^*_2)}{\partial \theta} =0
~$
and
$ \widetilde{A}_n \overset{a.s}{\underset{n,  \delta \rightarrow \infty}{\longrightarrow}} \E\big(\frac{\partial^2 \ell_{1,2}(\theta^*_2)}{\partial \theta^2} \big)$.

 \medskip
Since $ \E(\frac{\partial^2 \ell_{1,2}(\theta^*_2)}{\partial \theta^2})$  is a nonsingular matrix (see \cite{Francq2016}), we deduce from (\ref{mvt1.1}) that
\begin{equation}\label{mvt1.2}
\frac{n-t}{t-t^*}  (\widehat \theta_n(T_2) -\widehat \theta_n(T^*_2) ) \overset{a.s}{\underset{n,  \delta \rightarrow \infty}{\longrightarrow}} 0.
\end{equation}
 We conclude by the Taylor expansion on $\widehat L_n$ that
 \[
 \frac{1}{t-t^*}  | \widehat{L}_n (T_2,\widehat \theta_n(T_2) )  - \widehat{L}_n ( T_2,\widehat \theta_n(T^*_2) ) |
   \leq
   \frac{1}{2(t-t^*)} \| \widehat \theta_n(T_2)-\widehat \theta_n(T^*_2) \|^2 \sup_{\theta \in \Theta} \Big\|  \frac{\partial^2 \widehat L_n(T_2,\theta)}{\partial \theta^2} \Big\| \rightarrow 0~~~a.s.
   \]
   $~~~~~~~~~~~~~~~~~~~~~~~~~~~~~~~~~~~~~~~~~~~~~~~~~~~~~~~~~~~~~~~~~~~~~~~~~~~~~~~~~~~~~~~~~~~~~~~~~~~~~~~~~~~~~~~~~~~~~~~~ \blacksquare$
\end{enumerate}

\emph{\bf  Proof of Lemma  \ref{lem3}}\\
We detail the proof for the first order derivation ; the proof for the second order derivation follows the same reasoning.

 \medskip
Let $j \in \{ 1, \cdots,K^* \}$ and $l \in  \{1, \cdots , d \}$, we will show that
\begin{equation}\label{proof_lem3_partial}
\frac{1}{\sqrt{n^*_j}}\Big\|  \frac{\partial \widehat L_n(T^*_j,\theta)}{\partial \theta_l} - \frac{\partial \widehat L_{n,j}(T^*_j,\theta)}{\partial \theta_l}\Big\|_\Theta \overset{a.s}{\underset{n \rightarrow \infty}{\longrightarrow}} 0.
\end{equation}
Remark that
\begin{align*}
\frac{1}{\sqrt{n^*_j}}\Big\|  \frac{\partial \widehat L_n(T^*_j,\theta)}{\partial \theta_l} - \frac{\partial \widehat L_{n,j}(T^*_j,\theta)}{\partial \theta_l}\Big\|_\Theta
& \leq
 \frac{1}{\sqrt{n^*_j}} \sum_{t \in T^*_j}^{ } \Big\|  \frac{\partial \widehat \ell_t(\theta)}{\partial \theta_l} - \frac{\partial \widehat \ell_{t,j}(\theta)}{\partial \theta_l}\Big\|_\Theta\\
&\leq
    \frac{1}{\sqrt{n^*_j} }\sum_{t \in T^*_j}  \Big[Y_{t,j} \Big\|\frac{1}{\widehat{f}^\theta_t} \frac{\partial \widehat{f}^\theta_t}{\partial \theta_l}-\frac{1}{\widehat{f}^\theta_{t,j}} \frac{\partial \widehat{f}^\theta_{t,j}}{\partial \theta_l}\Big\|_\Theta +\Big\|\frac{\partial \widehat{f}^\theta_t}{\partial \theta_l}- \frac{\partial \widehat{f}^\theta_{t,j}}{\partial \theta_l}\Big\|_\Theta  \Big].
\end{align*}
Thus, by using the inequality $|a_1 b_1-a_2 b_2| \leq |a_1| | b_1-b_2| + |b_2| | a_1-a_2| $ $\forall a_1, a_2, b_1, b_2 \in \R$, we have
\begin{align*}
\frac{1}{\sqrt{n^*_j}}\Big\|  \frac{\partial \widehat L_n(T^*_j,\theta)}{\partial \theta_l} - \frac{\partial \widehat L_{n,j}(T^*_j,\theta)}{\partial \theta_l}\Big\|_\Theta
& \leq
\frac{1}{\sqrt{n^*_j} }\sum_{t \in T^*_j}  \Big[
Y_{t,j} \Big(
\Big\|\frac{1}{\widehat{f}^\theta_t} \Big\|_\Theta
  \Big\| \frac{\partial \widehat{f}^\theta_t}{\partial \theta_l}- \frac{\partial \widehat{f}^\theta_{t,j}}{\partial \theta_l}\Big\|_\Theta  + \Big\| \frac{\partial \widehat{f}^\theta_{t,j}}{\partial \theta_l} \Big\|_\Theta
  \Big\| \frac{1}{\widehat{f}^\theta_t} - \frac{1}{\widehat{f}^\theta_{t,j}}\Big\|_\Theta
  \Big)\\
  &~~~~~~~~~~~~~~~~~~~~~~+
\Big\|\frac{\partial \widehat{f}^\theta_t}{\partial \theta_l} -\frac{\partial \widehat{f}^\theta_{t,j}}{\partial \theta_l}\Big\|_\Theta
\Big]\\~\\
& \leq
\frac{1}{\sqrt{n^*_j} }\sum_{t \in T^*_j}  \Big[
\Big( \frac{Y_{t,j}}{\underline c}+1\Big)\Big\|\frac{\partial \widehat{f}^\theta_t}{\partial \theta_l} -\frac{\partial \widehat{f}^\theta_{t,j}}{\partial \theta_l}\Big\|_\Theta
+ \frac{1}{\underline c^2} Y_{t,j} \Big\| \frac{\partial \widehat{f}^\theta_{t,j}}{\partial \theta_l} \Big\|_\Theta \| \widehat{f}^\theta_t - \widehat{f}^\theta_{t,j} \|_\Theta
\Big] \\
& \leq
 C \frac{1}{\sqrt{n^*_j} }\sum_{t \in T^*_j}  \Big[
\Big(  Y_{t,j} +1   \Big)\Big\|\frac{\partial \widehat{f}^\theta_t}{\partial \theta_l} -\frac{\partial \widehat{f}^\theta_{t,j}}{\partial \theta_l}\Big\|_\Theta
+   Y_{t,j} \Big\| \frac{\partial \widehat{f}^\theta_{t,j}}{\partial \theta_l} \Big\|_\Theta  \| \widehat{f}^\theta_t - \widehat{f}^\theta_{t,j} \|_\Theta
\Big]
.
\end{align*}
 For $ t \in T^*_j$, from \textbf{A$_1 (\Theta)$}, we have
\begin{align}
  \Big\|  \frac{\partial \widehat{f}^\theta_t}{\partial \theta_l} - \frac{\partial \widehat{f}^\theta_{t,j}}{\partial \theta_l} \Big\|_\Theta
   \nonumber   &= \Big\| \frac{ \partial }{\partial \theta_l}  f( Y_{t-1},\cdots,Y_{1},0,\cdots;\theta )  - \frac{\partial   }{\partial \theta_l}  f( Y_{t-1,j},\cdots,Y_{t^*_{j-1}+1,j},0,\cdots ;\theta ) \Big\|_\Theta \\
     \label{proof_lem3_partial_f} & \leq \sum\limits_{k=t-t^*_{j-1}}^{t-1} \alpha^{(1)}_{k} |Y_{t-k}| \leq  \sum\limits_{k \geq t-t^*_{j-1} } \alpha^{(1)}_{k} Y_{t-k} .
\end{align}
Hence, from (\ref{proof_lem0_eq1}) and (\ref{proof_lem3_partial_f}), we get
\begin{multline*}
\frac{1}{\sqrt{n^*_j}}\Big\|  \frac{\partial \widehat L_n(T^*_j,\theta)}{\partial \theta_l} - \frac{\partial \widehat L_{n,j}(T^*_j,\theta)}{\partial \theta_l}\Big\|_\Theta
 \leq
 C \frac{1}{\sqrt{n^*_j} }\sum_{t \in T^*_j}
 \Big[
\Big(  Y_{t,j} +1 \Big) \sum\limits_{k \geq t-t^*_{j-1} } \alpha^{(1)}_{k} Y_{t-k}
+   Y_{t,j} \Big\| \frac{\partial \widehat{f}^\theta_{t,j}}{\partial \theta_l} \Big\|_\Theta   \sum\limits_{k \geq t-t^*_{j-1} } \alpha^{(0)}_{k} Y_{t-k}
\Big] . \\
  \leq
 C \frac{1}{\sqrt{n^*_j} }\sum_{\ell=1}^{n^*_j}  \Big[
\Big(  Y_{\ell + t^*_{j-1},j} +1 \Big) \sum \limits_{k \geq \ell } \alpha^{(1)}_{k} Y_{\ell + t^*_{j-1}-k}
+   Y_{\ell + t^*_{j-1},j} \Big\| \frac{\partial \widehat{f}^\theta_{\ell + t^*_{j-1},j}}{\partial \theta_l} \Big\|_\Theta \sum\limits_{k \geq \ell } \alpha^{(0)}_{k} Y_{\ell + t^*_{j-1}-k}
\Big] .
\end{multline*}
According to \cite{Kounias1969}, (\ref{proof_lem3_partial}) holds if
\begin{equation}\label{proof_lem3_Kounias_cond}
 \sum_{\ell \geq 1}  \frac{1}{\sqrt{\ell}} \E \Big[
\Big(  Y_{\ell + t^*_{j-1},j} +1 \Big) \sum \limits_{k \geq \ell } \alpha^{(1)}_{k} Y_{\ell + t^*_{j-1}-k}
+   Y_{\ell + t^*_{j-1},j} \Big\| \frac{\partial \widehat{f}^\theta_{\ell + t^*_{j-1},j}}{\partial \theta_l} \Big\|_\Theta \sum\limits_{k \geq \ell } \alpha^{(0)}_{k} Y_{\ell + t^*_{j-1}-k}
\Big] < \infty .
\end{equation}
 We have
 \[  \sum_{\ell \geq 1}  \frac{1}{\sqrt{\ell}} \E \Big[ \Big(  Y_{\ell + t^*_{j-1},j} +1 \Big) \sum \limits_{k \geq \ell } \alpha^{(1)}_{k} Y_{\ell + t^*_{j-1}-k} \Big]
 = C  \sum_{\ell \geq 1}  \frac{1}{\sqrt{\ell}}  \sum \limits_{k \geq \ell } \alpha^{(1)}_{k} < \infty .      \]
 Moreover,
\begin{align*}
 \Big\| \frac{\partial \widehat{f}^\theta_{\ell + t^*_{j-1},j}}{\partial \theta_l} \Big\|_\Theta
  &= \Big\| \frac{\partial  }{\partial \theta_l} f(0,\cdots ; \theta)  \Big\|_\Theta
  + \Big\| \frac{\partial  }{\partial \theta_l} f(Y_{\ell + t^*_{j-1} -1,j},\cdots Y_{t^*_{j-1} +1 , j},0,\cdots ; \theta)  - \frac{\partial  }{\partial \theta_l} f(0,\cdots ; \theta)  \Big\|_\Theta \\
  &\leq C + \sum_{k \geq 1} \alpha^{(1)}_{k} Y_{\ell + t^*_{j-1} -k,j}.
 \end{align*}
 Hence, from Minkowski's  inequality,
 \[ \Big\| \Big\| \frac{\partial \widehat{f}^\theta_{\ell + t^*_{j-1},j}}{\partial \theta_l} \Big\|_\Theta \Big\|_3
 = \Big( \E\Big[ \Big\| \frac{\partial \widehat{f}^\theta_{\ell + t^*_{j-1},j}}{\partial \theta_l} \Big\|_\Theta^3 \Big] \Big)^{1/3}
 \leq  C + \sum_{k \geq 1} \alpha^{(1)}_{k} \| Y_{\ell + t^*_{j-1} -k,j} \|_3 <  C(1 + \sum_{k \geq 1} \alpha^{(1)}_{k}) < C <\infty .  \]
 Thus, by using  Hölder and Minkowski's inequalities, it comes that
\begin{align*}
  \sum_{\ell \geq 1}  \frac{1}{\sqrt{\ell}} \E \Big[   Y_{\ell + t^*_{j-1},j} \Big\| \frac{\partial \widehat{f}^\theta_{\ell + t^*_{j-1},j}}{\partial \theta_l} \Big\|_\Theta \sum\limits_{k \geq \ell } \alpha^{(0)}_{k} Y_{\ell + t^*_{j-1}-k} \Big]
&\leq \sum_{\ell \geq 1}  \frac{1}{ \sqrt{\ell} }   \|Y_{\ell + t^*_{j-1},j} \|_3  \Big\| \Big\| \frac{\partial \widehat{f}^\theta_{\ell + t^*_{j-1},j}}{\partial \theta_l} \Big\|_\Theta  \Big\|_3
      \big\| \sum\limits_{k \geq \ell } \alpha^{(0)}_{k} Y_{\ell + t^*_{j-1}-k} \big\|_3 \\
&\leq C \sum_{\ell \geq 1}  \frac{1}{\sqrt{\ell}} \big\| \sum\limits_{k \geq \ell } \alpha^{(0)}_{k} Y_{\ell + t^*_{j-1}-k} \big\|_3 \\
&\leq C \sum_{\ell \geq 1}  \frac{1}{\sqrt{\ell}}  \sum\limits_{k \geq \ell } \alpha^{(0)}_{k} \|Y_{\ell + t^*_{j-1}-k} \|_3
\leq C \sum_{\ell \geq 1}  \frac{1}{\sqrt{\ell}}  \sum\limits_{k \geq \ell } \alpha^{(0)}_{k} < \infty .
\end{align*}
Hence, (\ref{proof_lem3_Kounias_cond}) is satisfied, and it holds that,
\begin{equation*}
\frac{1}{\sqrt{n^*_j}}\Big\|  \frac{\partial \widehat L_n(T^*_j,\theta)}{\partial \theta_l} - \frac{\partial \widehat L_{n,j}(T^*_j,\theta)}{\partial \theta_l}\Big\|_\Theta \overset{a.s}{\underset{n \rightarrow \infty}{\longrightarrow}} 0,\, \,\textrm{for~ any} \, \, j = 1, \cdots,K^* \, \, \textrm{and}\, \, l  \in \{1, \cdots , d \}.
\end{equation*}
Thus, (\ref{lem3_hatL_hatLj}) follows.

 \medskip
Remark that, one can go along the same lines as above by replacing $\widehat L_n$ by $L_n$ or $\widehat L_{n,j} $ by $L_{n,j} $; and obtain
\begin{equation}\label{proof_lem3_L_hatLj}
\frac{1}{\sqrt{n^*_j}} \Big\|  \frac{\partial L_n(T^*_j,\theta)}{\partial \theta} - \frac{\partial \widehat L_{n,j}(T^*_j,\theta)}{\partial \theta} \Big\|_\Theta \overset{a.s}{\underset{n \rightarrow \infty}{\longrightarrow}} 0,\, \,\textrm{for~ any} \, \, j =1, \cdots,K^* .
\end{equation}
\begin{equation}\label{proof_lem3_hatL_Lj}
\frac{1}{\sqrt{n^*_j}} \Big\|  \frac{\partial \widehat L_n(T^*_j,\theta)}{\partial \theta} - \frac{\partial L_{n,j}(T^*_j,\theta)}{\partial \theta} \Big\|_\Theta \overset{a.s}{\underset{n \rightarrow \infty}{\longrightarrow}} 0,\, \,\textrm{for~ any} \, \, j = 1, \cdots,K^* .
\end{equation}
Moreover, similar arguments as above can be easily applied to get
\[
\forall t \in \Z, ~ ~ \E \big\| \frac{\partial \ell_{t,j}(\theta) }{\partial \theta} \big\|_{\Theta} = C < \infty  ~ ~ ~ ~ \text{ for any } j=1,\cdots,K^* .
\]
By the uniform strong law of large number applied on the process  $\{  \frac{\partial \ell_{t,j}(\theta) }{\partial \theta} ,t\in \Z\}$, it holds that, for $ j=1,\cdots,K^* $,
\begin{equation}\label{proof_lem3_Lj_catLj}
\Big\| \frac{1}{n^*_j} \frac{\partial  L_{n,j}(T^*_j, \theta) }{\partial \theta}  - \frac{\partial  L^{(j)} (\theta) }{\partial \theta} \Big\|_\Theta
= \Big\| \frac{1}{n^*_j} \sum_{t\in T^*_j} \frac{\partial \ell_{t,j}(\theta) }{\partial \theta}  - \E (\frac{\partial \ell_{1,j}(\theta) }{\partial \theta} ) \Big\|_\Theta  \limitepsn 0 .
 \end{equation}
Thus, (\ref{lem3_hatL_calLj}) follows from (\ref{proof_lem3_hatL_Lj}) and (\ref{proof_lem3_Lj_catLj}).
This achieves the proof of the lemma.

$~~~~~~~~~~~~~~~~~~~~~~~~~~~~~~~~~~~~~~~~~~~~~~~~~~~~~~~~~~~~~~~~~~~~~~~~~~~~~~~~~~~~~~~~~~~~~~~~~~~~~~~~~~~~~~~~~~~~~~~~~~~~~~~~~~~~~~~~~~~~~~~~ \blacksquare$\\


\subsection{Proof of Theorem \ref{th3}}

Remark that for any $j \in \{1,\cdots,K^*\}$ $(\widehat{\theta}_n(\widehat{T}_j)-\theta^*_j)= (\widehat{\theta}_n(\widehat{T}_j)-\widehat{\theta}_n(T^*_j)) + (\widehat{\theta}_n(T^*_j)-\theta^*_j)$.
 According to Theorem \ref{th2}, it comes $(\widehat{t}_j - \widehat{t}^*_j) =o_P(\log (n))$.
 By relation (\ref{mvt1.2}), we obtain $(\widehat{\theta}_n(\widehat{T}_j)-\widehat{\theta}_n(T^*_j))=o_P(\frac{\log(n)}{n})$.    Therefore, $\sqrt{n^*_j} (\widehat{\theta}_n(\widehat{T}_j)-\widehat{\theta}_n(T^*_j)) \limiteproban 0 $ and it suffices to show that $\sqrt{n^*_j} (\widehat{\theta}_n(T^*_j)-\theta^*_j) \limiteloin \mathcal{N}_d(0,\Sigma_j)$ to conclude.\\

 Recall that for $n$ large enough, $\widehat{\theta}_n(T^*_j) \in \overset{\circ}{\Theta}$. By the mean value theorem, there exists $(\widetilde \theta_{n,k})_{1 \leq k \leq d} \in [\widehat{\theta}_n(T^*_j) , \theta^*_j]$ such that
 \begin{equation}\label{mvt2.1}
 \frac{\partial L_n (T^*_j, \widehat{\theta}_n(T^*_j))}{\partial \theta_k}= \frac{\partial L_n (T^*_j,\theta^*_j)}{\partial \theta_k} + \frac{\partial^2 L_n (T^*_j,\widetilde \theta_{n,k})}{\partial \theta \partial \theta_k} \left(\widehat{\theta}_n(T^*_j)-\theta^*_j\right).
\end{equation}
Let $J_n= -\Big(\frac{1}{n^*_j}\frac{\partial^2 L_n (T^*_j,\widetilde \theta_{n,k})}{\partial \theta \partial \theta_k}\Big)_{1 \leq k \leq d}$.
By relation (\ref{proof_lem3_L_hatLj}), Lemma \ref{lem3} and the proof of Theorem 2.2 of \cite{Francq2016}, we obtain $J_n \limitepsn J_j (\theta^*_j)$,
 where $J_j$ is the matrix defined in Theorem \ref{th3}.
 But, by Assumption (\textbf{A4}), it is easy to see that $J_j(\theta^*_j)$ is a non singular matrix. Thus, for $n$ large enough, $J_n$ is invertible and (\ref{mvt2.1}) gives

\begin{equation}\label{proof_th3_theta_L}
\sqrt{n^*_j} (\widehat{\theta}_n(T^*_j)-\theta^*_j)= -J^{-1}_n \Big[ \frac{1}{\sqrt{n^*_j}} \Big(  \frac{\partial L_n (T^*_j, \widehat{\theta}_n(T^*_j))}{\partial \theta} - \frac{\partial L_n (T^*_j,\theta^*_j)}{\partial \theta} \Big) \Big].
\end{equation}

From the proof of Theorem 2.2 of \cite{Francq2016}
, we get
\begin{equation}\label{proof_th3_Lj_Nd}
\frac{1}{\sqrt{n^*_j}} \frac{\partial L_{n,j} (T^*_j,\theta^*_j)}{\partial \theta} \limiteloin \mathcal{N}_d(0,I_j(\theta^*_j))
\end{equation}
where $I_j$ is given in Theorem  \ref{th3}.
According to (\ref{lem3_hatL_hatLj}), (\ref{proof_lem3_L_hatLj}) and (\ref{proof_lem3_hatL_Lj}), we have
\begin{equation*}
\frac{1}{\sqrt{n^*_j}}  \Big| \frac{\partial L_{n} (T^*_j,\theta^*_j)}{\partial \theta} - \frac{\partial L_{n,j} (T^*_j,\theta^*_j)}{\partial \theta}  \Big| \limitepsn 0 .
\end{equation*}
Hence,  (\ref{proof_th3_Lj_Nd}) implies
\begin{equation}\label{proof_th3_L_Nd}
 \frac{1}{\sqrt{n^*_j}} \frac{\partial L_{n} (T^*_j,\theta^*_j)}{\partial \theta} \limiteloin \mathcal{N}_d(0,I_j(\theta^*_j)).
\end{equation}
Moreover, from (\ref{lem3_hatL_hatLj}) and (\ref{proof_lem3_L_hatLj}), we have
\begin{equation*} 
\frac{1}{\sqrt{n^*_j}} \Big\| \frac{\partial L_{n} (T^*_j,\theta)}{\partial \theta} - \frac{\partial \widehat{L}_n (T^*_j,\theta)}{\partial \theta} \Big \|_{\Theta} \limitepsn 0 .
\end{equation*}
Therefore, since $ \partial \widehat L_n (T^*_j, \widehat{\theta}_n(T^*_j))/\partial \theta=0$, it follows that
\[
\frac{1}{\sqrt{n^*_j}} \frac{\partial L_n (T^*_j, \widehat{\theta}_n(T^*_j))}{\partial \theta}
= \frac{1}{\sqrt{n^*_j}}\Big(  \frac{\partial L_n (T^*_j, \widehat{\theta}_n(T^*_j))}{\partial \theta} - \frac{\partial \widehat L_n (T^*_j, \widehat{\theta}_n(T^*_j))}{\partial \theta} \Big) \limitepsn 0.
\]
Thus, uses (\ref{proof_th3_theta_L}) and (\ref{proof_th3_L_Nd}) to conclude.

$~~~~~~~~~~~~~~~~~~~~~~~~~~~~~~~~~~~~~~~~~~~~~~~~~~~~~~~~~~~~~~~~~~~~~~~~~~~~~~~~~~~~~~~~~~~~~~~~~~~~~~~~~~~~~~~~~~~~~~~~~~~~~~~~~~~~~~~~~~~~~~~~ \blacksquare$\\


\begin{thebibliography}{99}

\bibitem{Francq2016}
{\sc Ahmad, A. and Francq, C.}
\newblock Poisson QMLE of count time series models.
\newblock {\em Journal of Time Series Analysis 37}, (2016),  291-314.

\bibitem{Arlot2016}
{\sc Arlot, S., Celisse, A. and Celisse, Z.}
\newblock A kernel multiple change-point algorithm via model selection.
\newblock {\em arXiv:1202.3878v2 : 2016}.


\bibitem{Arlot2009}
{\sc Arlot, S. and Massart, P.}
\newblock Data-driven Calibration of Penalties for Least-Squares Regression.
\newblock {\em Journal of Machine Learning Research 10}, (2009),  245-279.

\bibitem{Bai1997}
{\sc Bai, J.}
\newblock Estimating multiple breaks one at a time.
\newblock {\em Econometric theory 13 (3)}, (1997), 315-352.


\bibitem{Bai1998}
{\sc Bai, J. and Perron, P.}
\newblock Estimating and testing linear models with multiple structural changes.
\newblock {\em Econometrica 66}, (1998),  47-78.

\bibitem{Bardet2012}
{\sc Bardet, J. M., Kengne, K. and Wintenberger, O.}
\newblock Multiple breaks detection in general causal time series using penalized quasi-likelihood.
\newblock {\em Electronic Journal of Statistics 6}, (2012), 435-477.


\bibitem{Baudry2010}
{\sc Baudry, J. P., and Maugis, C. and Michel, B.}
\newblock Slope Heuristics: overview and implementation.
\newblock {\em RR-INRIA $n^o$7223}, 2010.











\bibitem{Cleynen2014}
{\sc  Cleynen, A. and Lebarbier, E.}
\newblock Segmentation of the Poisson and negative binomial rate models: a penalized estimator.
\newblock {\em ESAIM: Probability and Statistics 18}, (2014), 750-769.

\bibitem{Cleynen2017}
{\sc  Cleynen, A. and Lebarbier, E.}
\newblock Model selection for the segmentation of multiparameter exponential family distributions.
\newblock {\em Electronic Journal of Statistics 11}, (2017), 800-842.







\bibitem{Davis2016}
{\sc Davis, R. A., Hancock, S. A. and Yao, Y.-C.}
\newblock On consistency of minimum description length model selection for piecewise autoregressions.
\newblock {\em Journal of Econometrics 194}, (2016) 360-368.

\bibitem{David2008}
 \textsc{Davis, R. A., Lee, T. C. M. and Rodriguez-Yam, G. A.}
\newblock  Break detection for a class of nonlinear time  series models.
\newblock \textit{Journal of Time Series Analysis 29 (5)},  (2008), 834-867.


\bibitem{Davis2013}
{\sc Davis, R. A. and Yau, C. Y.}
\newblock Consistency of minimum description length model selection for piecewise stationary time series models.
\newblock {\em Electronic Journal of Statistics 7}, (2013) 381-411.


\bibitem{Diop2017}
{\sc Diop, M.L. and Kengne, W.}
\newblock Testing parameter change in general integer-valued time series.
\newblock {\em J. Time Ser. Anal. 38}, (2017) 880-894.


\bibitem{Douc2017}
{\sc Douc, R., Fokianos, K., and Moulines, E.}
\newblock Asymptotic properties of quasi-maximum likelihood estimators in observation-driven time series models.
\newblock {\em Electronic Journal of Statistics  11}, (2017),  2707-2740.



 \bibitem{Doukhan2015}
{\sc Doukhan, P. and Kengne, W. }
\newblock Inference and testing for structural change in general poisson autoregressive models.
\newblock {\em Electronic Journal of Statistics 9}, (2015), 1267-1314.




 \bibitem{Fokianos2009}
 {\sc Fokianos, K., Rahbek, A. and Tjøstheim, D.}
 \newblock Poisson autoregression.
 \newblock {\em Journal of the American Statistical Association 104}, (2009), 1430-1439.

 \bibitem{Fokianos2011}
{\sc Fokianos, K. and Tjøstheim, D.}
\newblock Log-linear Poisson autoregression.
\newblock {\em Journal of multivariate analysis 102},  (2011) 563-578.

 \bibitem{Fokianos2012}
 {\sc Fokianos, K., and Tjøstheim, D.}
 \newblock Nonlinear Poisson autoregression.
 \newblock {\em Ann. Inst. Stat. Math. 64}, (2012) 1205-1225.



  \bibitem{Fryzlewicz2014}
 {\sc Fryzlewicz, P. }
 \newblock Wild binary segmentation for multiple change-point detection.
 \newblock {\em The Annals of Statistics 42 (6)}, (2014) 2243-228.


  \bibitem{Fryzlewicz2014b}
 {\sc Fryzlewicz, P., and Subba Rao, S.}
 \newblock Multiple-change-point detection for auto-regressive conditional heteroscedastic processes.
 \newblock {\em Journal of the Royal Statistical Society: series B 76 (5)}, (2014) 903-924.


 \bibitem{Harchaoui2010}
 {\sc Harchaoui, Z. and L\'evy-Leduc, C. }
 \newblock Multiple Change-Point Estimation With a Total Variation Penalty.
 \newblock {\em Journal of the American Statistical Association  105}, (2010) 1480-1493.


\bibitem{Hudecova2013}
 {\sc Hudecov{\'a}, {\v{S}}.}
 \newblock Structural changes in autoregressive models for binary time series.
 \newblock {\em Journal of Statistical Planning and Inference 143}, (2013) 1744-1752.


 \bibitem{Inclan1994}
 {\sc Inclán, C. and Tiao, G. C.}
 \newblock  Use of cumulative sums of squares for retrospective detection of changes of variance.
 \newblock {\em Journal of the American Statistical Association 89 (427)}, (1994), 913-923.







\bibitem{Kounias1969}
{\sc Kounias, E.G. and Weng, T.-S }
\newblock  An inequality and almost sure convergence.
\newblock {\em Annals of Mathematical Statistics 33}, (1969), 1091-1093.




\bibitem{Lavielle2000b}
{\sc  Lavielle, M. and Ludena, C. }
\newblock The multiple change-points problem for the spectral distribution.
\newblock {\em Bernoulli 6}, (2000), 845-869.


\bibitem{Lebarbier2005}
{\sc Lebarbier E. }
\newblock  Detecting multiple change-points in the mean of Gaussian process by model selection.
\newblock {\em Signal Processing 85}, (2005), 717-736.



\bibitem{Lerasle2011}
{\sc Lerasle, M. }
\newblock   Optimal model selection for density estimation of stationary data under various mixing conditions.
\newblock {\em  The Annals of Statistics 39 (4)}, (2003), 1852-1877.









\bibitem{Yau2016}
{\sc Yau, C.Y. and Zhao, Z.}
 \newblock Inference for multiple change points in time series via likelihood ratio scan statistics.
 \newblock {\em Journal of the Royal Statistical Society: Series B 78}, (2016), 895-916.










\end{thebibliography}
\end{document}